\documentclass[a4paper,10pt]{article}

\usepackage{amsmath,amsthm}
\usepackage{upref}
\usepackage{amsfonts}
\usepackage{pstricks}
\usepackage{amssymb}

\newtheorem{lmm}{Lemma}[section]

\newtheorem{prp}{Proposition}[section]

\newtheorem{thm}{Theorem}[section]

\newcommand{\sequence}{\operatorname{sequence}}
\newcommand{\Irr}[1]{\operatorname{Irr}(#1)}

\def\[{[\![}
\def\]{]\!]}

\newtheorem{re}{Remark}[section]

\newtheorem{df}{Definition}[section]
\newtheorem{exa}{Example}[section]

\title{On the  parametrization of the simple modules for 
Ariki-Koike algebras at roots of unity}

\author{Nicolas Jacon\footnote{Institut Girard Desargues, Lyon I}}
\date{}

\begin{document}

\maketitle

\begin{abstract}                        
Following the ideas of M.Geck and R.Rouquier in \cite{geckrouq},
 we show that there exists a ``canonical basic set'' of Specht 
modules in bijection with the set of simple modules of Ariki-Koike 
algebras at roots of unity. Moreover, we determine the parametrization 
of this set and we give the consequences of these results on the 
representation theory of Ariki-Koike algebras. 
\end{abstract}

\section{Introduction}

 Let $H$ be an Iwahori-Hecke algebra of a finite Weyl group $W$ 
(or of an extended Weyl group) defined over $A:=\mathbb{Z}[v,v^{-1}]$ 
where $v$ is an indeterminate. For any ring homomorphism $\theta:A\to{k}$
 into a field $k$, we have the corresponding specialized algebra
 $H_k:=k\otimes_A{H}$. Assume that $K$ is the field of fractions  of $A$. 
Then, the algebra $H_K$ is split semi-simple and isomorphic to the group
 algebra $K[W]$.

One of the major problems is the determination of the simple $H_k$-modules
 when $H_k$ is not semi-simple. This problem may be attacked by using 
the corresponding decomposition matrix, which relates the simple 
$H_k$-modules to the simple $H_K$-modules, by a process of modular
 reduction.

In \cite{geckrouq}, using an ordering of simple $H_K$-modules by  
Lusztig $a$-function, Geck and Rouquier showed that there exists a
 canonical set $\mathcal{B}\subset{\Irr{H_K}}$ in bijection with 
$\Irr{H_k}$ (see also \cite{geck2} and  \cite{geck3}). As a consequence,
 this set gives a natural way for labeling the simple $H_k$-modules. 
The  Lusztig $a$-function has several interpretations: one in terms 
of  Kazhdan-Lusztig basis of $H$ and one in terms of Schur elements
 (the coincidence is shown in \cite{lus2}).

In this paper, we consider the case of Ariki-Koike algebras. 
Let $d\in{\mathbb{N}_{>0}}$ and $\mathcal{H}_{n}:=\mathcal{H}_n
(v;u_0,...,u_{d-1})$ be the Ariki-Koike algebra of type $G(d,1,n)$
 where $v$, $u_0$, $u_1$,..., $u_{d-1}$ are $d+1$ parameters. 
This algebra appeared independently in \cite{AriKoi}    and in 
\cite{BrouMal} and can be seen as an analogue of the Iwahori-Hecke 
algebra for the complex reflection group  $(\mathbb{Z}/d\mathbb{Z})
\wr{\mathfrak{S}_n}$.

As the  Iwahori-Hecke algebras of type $A_{n-1}$ and $B_n$ are 
special cases of Ariki-Koike algebras, it is a natural question
 to ask whether the set $\mathcal{B}$ is well-defined for the 
 Ariki-Koike algebras. We don't have Kazhdan-Lusztig type basis
 for  Ariki-Koike algebras but we do have Schur elements. Hence, 
we can also define an $a$-function for this type of algebras.

When the parameters  are generic, the representation theory of 
 $\mathcal{H}_{n}$ over a field has been studied in \cite{AriKoi}.
 It was shown that this algebra is semi-simple and that the 
simple modules $S^{\underline{\lambda}}$ are parametrized by
 the $d$-tuples of partitions $\underline{\lambda}=
(\lambda^{(0)},...,\lambda^{(d-1)})$ of rank $n$.

Now, we consider the following choice of parameters:
$$v=\eta_e,\qquad{u_i =\eta_e^{v_j},\ \ j=0,...,d-1,}$$
where $\eta_e:=\textrm{exp}{(\frac{2i\pi}{e})}\in\mathbb{C}$ 
 and where $0\leq{v_0}\leq ...\leq v_{d-1}<e$.

In this case,  $\mathcal{H}_{n}$ is  not semi-simple in general
 (it may be semi-simple if $n$ is small). We can also attach to 
each $d$-partition $\underline{\lambda}$ an $\mathcal{H}_n$-module
 $S^{\underline{\lambda}}$ but it is  reducible in general. 
 In \cite{grahale}, Graham and Lehrer constructed a bilinear 
form on each $S^{\underline{\lambda}}$ and proved that the 
simple modules of $\mathcal{H}_{n}$ are given by the non 
zero $D^{\underline{\lambda}}:=S^{\underline{\lambda}}/
\textrm{rad}(S^{\underline{\lambda}})$. Then, Ariki and 
Mathas (see \cite{Ari4} and \cite{AriMat1}) described the
 non zero   $D^{\underline{\lambda}}$ and showed that they 
are indexed by some ``Kleshchev'' multipartitions.
Unfortunately, we only know a recursive description of 
this kind of multipartitions.

The aim of this paper is to prove that there exists a set 
$\mathcal{B}$ labeled by simple modules of a semi-simple 
Ariki-Koike algebra which satisfy the same property as the 
canonical set of \cite{geckrouq} for Iwahori-Hecke algebras. 
Moreover, we prove that the parametrization of this set 
coincides with the parameterization of the set of simple 
$\mathcal{H}_{n}$-modules found by  Foda, Leclerc, Okado, 
Thibon and Welsh  (see \cite{FLOTW}).  
The proof requires results of Ariki and Foda et al. about 
Ariki-Koike algebras,  some properties of $a$-function and 
some combinatorial objects such as symbols.

 These results have several consequences. First, it extends 
some results proved in \cite{geck2} where the set $\mathcal{B}$
 was determined for Hecke algebras of type $A_{n-1}$ and in 
   \cite{a-fonction} where  $\mathcal{B}$ was determined for
 Hecke algebras of type $D_n$ and $B_n$ with the following 
diagram:
\\
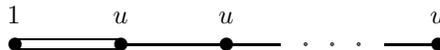
\begin{center}
\begin{picture}(240,20)
\put( 50,10){\circle*{5}}
\put( 47,18){$1$}
\put( 50,8){\line(1,0){40}}
\put( 50,12){\line(1,0){40}}
\put( 90,10){\circle*{5}}
\put( 87,18){$u$}
\put( 90,10){\line(1,0){40}}
\put(130,10){\circle*{5}}
\put(127,18){$u$}
\put(130,10){\line(1,0){20}}
\put(160,10){\circle{1}}
\put(170,10){\circle{1}}
\put(180,10){\circle{1}}
\put(190,10){\line(1,0){20}}
\put(210,10){\circle*{5}}
\put(207,18){$u$}
\end{picture}
\end{center}

 In particular, it yields  the determination of the set $\mathcal{B}$ 
for  Hecke algebras of type $B_n$ and, hence, completes the 
classification of  the canonical basic set for Hecke algebras
 of finite Weyl groups with one parameter.

Moreover, using the results of Ariki and the ideas developed by 
Lascoux, Leclerc and Thibon in \cite{LLT} for Hecke algebras of 
type $A_{n-1}$, we obtain a ``triangular'' algorithm for computing
 the decomposition matrix for Ariki-Koike algebras when the 
parameters are roots of unity.

The paper is divided in six parts. First, we briefly summarize 
the results of Ariki and Foda et al. about simple modules for 
Ariki-Koike algebras and links with quantum groups. Then, in 
the second part, following \cite{BrouKim},  we introduce Schur 
elements and $a$-function associated to the simple modules of  
a semi-simple  Ariki-Koike algebra. In the third part, we prove 
some combinatorial properties about the multipartitions of Foda 
et al. and the $a$-function. The fourth part contains the main 
theorem of the paper: we give an interpretation of the multipartitions
 of  Foda et al. in terms of the $a$-function. Finally, in the 
fifth  part, we study a particular case of Ariki-Koike algebras 
: the case of Hecke algebras of type $B_n$ with equal parameters.

\section{Preliminaries}

Let  $R$ be a commutative associative ring with unit and let 
 $v$, $u_0$, ..., $u_{d-1}$ be $d+1$ invertible elements in $R$.
 Let $n\in{\mathbb{N}}$. We define the  Ariki-Koike algebra 
$\mathcal{H}_{R,n}:=\mathcal{H}_{R,n}{({v};{u}_0,...,{u}_{d-1})}$
 over $R$ to be the unital associative $R$-algebra generated by
 $T_0$, $T_1$,..., $T_{n-1}$  subject to the relations:
\begin{align*}
 &   (T_i-{v})(T_i+1)=0& & \textrm{for}\    1\leq{i}\leq{n-1},&\\
 &(T_0-{u}_0)(T_0-{u}_1)...(T_0-{u}_{d-1})=0,&& &\\
 &  T_iT_{i+1}T_{i}=T_{i+1}T_iT_{i+1}& &\textrm{for}\       
 1\leq{i}\leq{n-2}, &  \\
 & T_iT_j=T_jT_i& &\textrm{for}\ |i-j|>1, &\\
  & T_0T_1T_0T_1=T_1T_0T_1T_0.&& & 
 \end{align*}
These relations are obtained by deforming the relations of the 
wreath product $(\mathbb{Z}/d\mathbb{Z})\wr{\mathfrak{S}_n}$. 
The last three ones  are known as type $B$ braid relations. 

 It is known that the simple modules of $(\mathbb{Z}/d\mathbb{Z})
\wr{\mathfrak{S}_n}$   are indexed by the $d$-tuples of partitions.
 We will see that the same is true for the semi-simple Ariki-Koike
 algebras defined over a field. In this paper, we say that  
$\underline{\lambda}$  is a $d$-partition of rank $n$ if:
\begin{itemize}
\item  $\underline{\lambda}=(\lambda^{(0)},..,\lambda^{(d-1)})$ 
where, for $i=0,...,d-1$, $\lambda^{(i)}=(\lambda^{(i)}_1,...,
\lambda^{(i)}_{r_i})$ is a  partition of rank $|\lambda^{(i)}|$
  such that $\lambda^{(i)}_1\geq{...}\geq{ \lambda^{(i)}_{r_i}}>0$,
\item $\displaystyle{\sum_{k=0}^{d-1}{|\lambda^{(k)}|}}=n$.
\end{itemize}
We denote by $\Pi_{n}^d$ the set of  $d$-partitions of rank $n$.

For each $d$-partition $\underline{\lambda}$ of rank $n$, we 
can associate a $\mathcal{H}_{R,n}$-module $S^{\underline{\lambda}}$ 
which is free over $R$. This is called  a Specht 
module\footnote{ Here, we use the definition of the 
classical Specht modules. Note that the  results in  
\cite{DJMATHAS} are given in terms of dual Specht modules. 
The passage from classical Specht modules to their duals is
 provided by the map $(\lambda^{(0)},\lambda^{(1)},...,\lambda^{(d-1)})
\mapsto{(\lambda^{(d-1)'},\lambda^{(d-2)'},...,\lambda^{(0)'})}$ 
where, for $i=0,...,d-1$,  $\lambda^{(i)'}$ denotes the conjugate 
partition.}.

Assume that $R$ is a field. 
 Then, for each $d$-partition of rank $n$, there is a natural 
bilinear form which is defined over each $S^{\underline{\lambda}}$.
 We denote by $\textrm{rad}$ the radical associated to this bilinear form. 
The non zero  $D^{\underline{\lambda}}:=S^{\underline{\lambda}}/
\textrm{rad}{(S^{\underline{\lambda}})}$ form a complete set of
 non-isomorphic simple $\mathcal{H}_{R,n}$-modules (see for example
 \cite[chapter 13]{Arilivre}). In particular, if $\mathcal{H}_{R,n}$
 is semi-simple, we have $\textrm{rad}(S^{\underline{\lambda}})=0$
 for all $\underline{\lambda}\in\Pi_{n}^d$ and  the set of simple 
modules are given by the $S^{\underline{\lambda}}$.

Now, we explain in more details the representation theory of  
Ariki-Koike algebras.   We summarize the results of Ariki and 
Foda et al. following \cite{AriMat1} and \cite{FLOTW}.

\subsection{Simple modules  of Ariki-Koike algebras and decomposition 
numbers}
Here, we assume that $R$ is a field and we consider the Ariki-Koike 
algebra $\mathcal{H}_{R,n}$ defined over $R$ with parameters $v$, 
$u_0$,..., $u_{d-1}$. First, we have a criterion of semi-simplicity:

\begin{thm}[Ariki \cite{Ari3}]\label{semisimple}  $\mathcal{H}_{R,n}$ 
is split semi-simple if and only if we have:
\begin{itemize}
\item for all $i\neq{j}$ and for all   $d\in{\mathbb{Z}}$ such 
that $|d|<{n}$, we have:
$$v^d u_i\neq{u_j},$$
\item $\displaystyle{\prod_{i=1}^{n}{(1+v+...+v^{i-1})}\neq{0}}.$
\end{itemize}
\end{thm}
In this case, as noted at the beginning of this section, Ariki 
proved that the simple modules  are labeled by the $d$-partitions
 $\underline{\lambda}$ of rank $n$.

Using the results of Dipper and Mathas (\cite{DipMat}),  the case 
where $\mathcal{H}_{R,n}$ is  not semi-simple  can be  reduced to 
the case where all the $u_i$ are powers of $v$. 

 Here, we assume that $R$ is a field of characteristic $0$ and 
that $v$ is a primitive $e^{\textrm{th}}$-root of  unity with 
$e\geq 2$: 
\begin{eqnarray*}
&&u_j=\eta_e^{v_{j}}\ \textrm{for}\ j=0,...,d-1,\\
&&v=\eta_{e}, 
\end{eqnarray*}
where $\eta_e:=\textrm{exp}(\frac{2i\pi}{e})$  and where 
$0\leq{v_0}\leq ...\leq v_{d-1}<e$. The problem is to find
 the non zero $D^{\underline{\lambda}}$. To solve it, it 
is convenient to use the language of decomposition matrices.

 We define $R_0 (\mathcal{H}_{R,n})$ to be  the Grothendieck
 group of finitely generated $\mathcal{H}_{R,n}$-modules, 
which is generated by the simple modules of $\mathcal{H}_{R,n}$.
 For  a   finitely generated $\mathcal{H}_{R,n}$-module $M$,
 we denote by $[M]$ its equivalence class in $R_0 (\mathcal{H}_{R,n})$.

Now, let $R_1 (\mathcal{H}_{R,n})$ be the  Grothendieck group 
of finitely generated projective $\mathcal{H}_{R,n}$-modules. 
 This is generated by the  indecomposable  projective 
$\mathcal{H}_{R,n}$-modules. For a finitely generated projective
 $\mathcal{H}_{R,n}$-module $P$, we denote by $[P]_p$ its 
equivalent class in $R_1 (\mathcal{H}_{R,n})$. We have an 
injective homomorphism
$$c : R_1 (\mathcal{H}_{R,n})\to{R_0 (\mathcal{H}_{R,n})}$$
given    by $c([P]_p)=[P]$ where $P$ is  a finitely generated 
projective $\mathcal{H}_{R,n}$-module.

Following \cite{AriMat1}, we denote by $\mathcal{F}_n$ the free
 abelian group with $\mathbb{C}$-basis the set 
$\{\[S^{\underline{\lambda}}\]\ |\ \underline{\lambda}\in{\Pi_n^d}\}$.
  $\mathcal{F}_n$  can be seen as the Grothendieck group of a 
semi-simple Ariki-Koike algebra.

 Let  $\Phi^d_n=\{ \underline{\mu}\in{\Pi^d_n} |\ D^{\underline{\mu}}
\neq{0}\}$.  $R_0(\mathcal{H}_{R,n})$ is generated by the set 
$\{[D^{\underline{\mu}}]|\ \underline{\mu}\in{\Phi^d_n}\}$. 
Hence, for all $\underline{\lambda}\in{\Pi_n^d}$, there exist
 numbers $d_{\underline{\lambda},\underline{\mu}}$ with 
$\underline{\mu}\in{\Phi^d_n}$ such that :
$$[S^{\underline{\lambda}}]=\sum_{\underline{\mu}\in{\Phi^d_n}}
{d_{\underline{\lambda},\underline{\mu}}[D^{\underline{\mu}}]}.$$
 The matrix $(d_{\underline{\lambda},\underline{\mu}})
_{\underline{\lambda}\in{\Pi_n^d},\underline{\mu}\in{\Phi^d_n}}$ 
is called the decomposition matrix of $\mathcal{H}_{R,n}$.

Hence, we have a homomorphism:
$$d:\mathcal{F}_n\to{ R_0 (\mathcal{H}_{R,n})}, $$
defined by $\displaystyle{d(\[S^{\underline{\lambda}}\])=
[S^{\underline{\lambda}}]=\sum_{\underline{\mu}\in{\Phi^d_n}}
{d_{\underline{\lambda},\underline{\mu}}[D^{\underline{\mu}}]}}$.
 $d$ is called the decomposition map. Note that if 
$\mathcal{H}_{R,n}$ is semi-simple, the decomposition matrix 
is just the identity.

 Now,  by Brauer reciprocity, we can see  that the  indecomposable 
 projective $\mathcal{H}_{R,n}$-modules are labeled by  $\Phi^d_n$
 and that there exists an injective homomorphism:
$$e: R_1 (\mathcal{H}_{R,n})\to{\mathcal{F}_n},$$
such that  for $\underline{\mu}\in{\Phi^d_n}$, if
 $P^{\underline{\mu}}$  is the   indecomposable  projective 
$\mathcal{H}_{R,n}$-module which is the  projective cover of 
$D^{\underline{\mu}}$, we have:
$$e([P^{\underline{\mu}}]_p)=\sum_{\underline{\lambda}\in{ \Pi^d_n}}
   d_{\underline{\lambda},\underline{\mu}}\[S^{\underline{\lambda}}\].$$
By \cite[Theorem 3.7]{grahale}, we obtain the following commutative
 diagram:\\
\begin{center}
\begin{picture}(200,60)
\put(60,40){\vector(-2,-1){30}}
\put(105,40){\vector(2,-1){30}}
\put(62,43){$R_1 (\mathcal{H}_{R,n})$}
\put(20,15){$\mathcal{F}_n$}
\put(33,18){\vector(1,0){100}}
\put(135,15){$R_0 (\mathcal{H}_{R,n})$}
\put(43,35){$e$}
\put(118,35){$c$}
\put(81,20){$d$}
\end{picture}
\end{center}

 Now, we turn to the problem of determining which 
 $D^{\underline{\lambda}}$ are non zero. Ariki and Mathas solve it 
  by using deep results about quantum groups:

\subsection{Links with quantum groups and Kleshchev
 $d$-partitions}\label{combinatoire}
We keep the notations of the first paragraph of this section.
First, we introduce some notations and we define  Kleshchev 
$d$-partitions following \cite{AriMat1}.

Let  $\underline\lambda={(\lambda^{(0)},...,\lambda^{(d-1)})}$
 be a $d$-partition of rank  $n$. The diagram of  
$\underline{\lambda}$ is the following set:
$$[\underline{\lambda}]=\left\{ (a,b,c)\ |\ 0\leq{c}\leq{d-1},\ 
1\leq{b}\leq{\lambda_a^{(c)}}\right\}.$$

The elements of this diagram are called   the nodes of  
$\underline{\lambda}$. Let  $\gamma=(a,b,c)$ be a node of  
$\underline{\lambda}$. The residue of  $\gamma$  associated 
to the set  $\{e;{v_0},...,{v_{d-1}}\}$ is the element of 
$\mathbb{Z}/e\mathbb{Z}$ defined by:
$$\textrm{res}{(\gamma)}=(b-a+v_{c})(\textrm{mod}\ e).$$

If $\gamma $ is a node with residue $i$, we say that  $\gamma$ 
is an  $i$-node. 
Let  $\underline{\lambda}$ and  $\underline{\mu}$ be two $d$-partitions
 of rank  $n$ and $n+1$ such that  $[\underline{\lambda}]\subset
{[\underline{\mu}]}$. There exists a node $\gamma$ such that  
$[\underline{\mu}]=[\underline{\lambda}]\cup{\{\gamma\}}$. Then,
 we denote $[\underline{\mu}]/[\underline{\lambda}]=\gamma$. If 
$\textrm{res}{(\gamma)}=i$, we say that  $\gamma$ is an addable 
$i$-node  for   $\underline{\lambda}$ and a removable  $i$-node 
for  $\underline{\mu}$.

Now, we consider  the following order on the set of removable and 
addable nodes of a $d$-partition:  
 we say that $\gamma=(a,b,c)$ is below  $\gamma'=(a',b',c')$ if  
$c<c'$ or if  $c=c'$ and $a<a'$.

This order will be called the AM-order and the notion of normal nodes
 and good nodes below are linked with this order (in the next paragraph, 
we will give another order on the set of nodes which is distinct from  
this one).

Let   $\underline{\lambda}$ be a   $d$-partition and  let $\gamma$ 
be  an $i$-node, we say that  $\gamma$ is a normal $i$-node of  
$\underline{\lambda}$ if, whenever $\eta$ is an $i$-node of 
$\underline{\lambda}$   below  $\gamma$, there are more removable
 $i$-nodes between $\eta$ and $\gamma$ than addable $i$-nodes 
between  $\eta$ and $\gamma$. 
If $\gamma$ is the highest normal $i$-node of    $\underline{\lambda}$, 
we say that  $\gamma$ is a good $i$-node.

We can now define the notion of Kleshchev $d$-partitions  associated 
to the set  $\{e;{v_0},...,{v_{d-1}}\}$:
\begin{df}
The Kleshchev $d$-partitions are defined   recursively  as follows.
\begin{itemize}
   \item The empty partition $\underline{\emptyset}:=(\emptyset,
\emptyset,...,\emptyset)$ is  Kleshchev.
    \item If $\underline{\lambda}$ is  Kleshchev, there exist 
$i\in{\{0,...,e-1\}}$ and a good $i$-node $\gamma$ such that if 
we remove $\gamma$ from  $\underline{\lambda}$, the resulting  
$d$-partition is Kleshchev.
\end{itemize}

We denote  by $\Lambda^0_{\{e;v_0,...,v_{l-1}\}}$ the set of  
Kleshchev $d$-partitions associated to  the set  $\{e;{v_0},.
..,{v_{d-1}}\}$. If there is no ambiguity concerning  $\{e;{v_0},
...,{v_{d-1}}\}$, we denote it by $\Lambda^0$.

\end{df}

 Now, let  $\mathfrak{h}$ be a free $\mathbb{Z}$-module with 
basis $\{h_i,\mathfrak{d}\ |\ 0\leq{i}<e\}$ and let $\{\Lambda_i,
\delta\ |\ 0\leq{i}<e\}$ be the dual basis  with respect to the
 pairing: 
$$\langle\ ,\ \rangle :\mathfrak{h}^*\times{ \mathfrak{h}}\to{\mathbb{Z}}$$
such that $\langle\Lambda_i,h_j\rangle=\delta_{ij}$, 
$\langle\delta,\mathfrak{d}\rangle=1$ and  $\langle\Lambda_i,
\mathfrak{d}\rangle=\langle\delta,h_j\rangle=0$ for $0\leq{i,j}<e$.
 For $0\leq{i}<e$, we define the simple roots of $\mathfrak{h}^*$ by:
$$\alpha_i=\left\{ \begin{array}{ll} 2\Lambda_0-\Lambda_{e-1}-
\Lambda_1+\delta & \textrm{if}\  i=0, \\
                                   2\Lambda_i-\Lambda_{i-1}-
\Lambda_{i+1} & \textrm{if}\  i>0, \\
    \end{array} \right.$$
where $\Lambda_e:=\Lambda_0$. The  $\Lambda_i$ are called  the 
fundamental weights.

Now, let $q$ be an indeterminate and let $\mathcal{U}_q$ be the 
quantum group of type $A^{(1)}_{e-1}$.
This is a unital associative algebra over  $\mathbb{C}(q)$ which 
is generated by elements $\{e_i,f_i\ |\ i\in{\{0,...,e-1\}}\}$ 
and $\{k_h\ |\ h\in{\mathfrak{h}}\}$ 
(see for example \cite[chapter 6]{mathas} for the relations).

For $j\in{\mathbb{N}}$ and $l\in{\mathbb{N}}$, we define:
\begin{itemize}
\item $\displaystyle{[j]_q:=\frac{q^j-q^{-j}}{q-q^{-1}}}$,
\item $\displaystyle{ [j]_q^{!}:=[1]_q[2]_q...[j]_q}$,
\item  $\displaystyle{   \left[ \begin{array}{c} l      \\ 
 j \end{array} \right]_q=\frac{[l]_q^!}{ [j]_q^![l-j]_q^!}}$.
\end{itemize}

We now introduce the Ariki's theorem. For details, we refer 
to \cite{Arilivre}. Let $\mathcal{A}=\mathbb{Z}[q,q^{-1}]$. 
We consider the Kostant-Lusztig   $\mathcal{A}$-form of 
$\mathcal{U}_q$  which is denoted by  $\mathcal{U}_{\mathcal{A}}$:
 this is a  $\mathcal{A}$-subalgebra of $\mathcal{U}_q$ generated 
by the divided powers  $e_i^{(l)}:=\displaystyle{\frac{e^{l}_i}{[l]^{!}_q}}$
, $f_j^{(l)}:=\displaystyle{\frac{f^{l}_j}{[l]^{!}_q} }$  for 
 $0\leq{i,j}<e$ and $l\in{\mathbb{N}}$ and by  $k_{h_{i}}$, 
$k_{\mathfrak{d}}$,  $k^{-1}_{h_i}$, $k^{-1}_{\mathfrak{d}}$ for
 $0\leq{i}<e$ (see \cite[\S 3.1]{lusztig}). Now, if $S$ is a ring 
and $u$   an invertible element in $S$, we can form the specialized
 algebra $\mathcal{U}_{S,u}:=S\otimes_{\mathcal{A}}\mathcal{U}_
{\mathcal{A}}$  by specializing the indeterminate $q$ to $u\in{S}$.

 Let $\underline{\lambda}$ and $\underline{\mu}$ be two $d$-partitions
 of rank $n$ and $n+1$ such that  there exists an $i$-node $\gamma$
 such that  $[\underline{\mu}]=[\underline{\lambda}]\cup{\{\gamma\}}$.
 We define:
\begin{align*}
N_i^{a}{(\underline{\lambda},\underline{\mu})}=&   \sharp\{ \textrm{addable }
\ i{}-\textrm{nodes of } \underline{\lambda}\ \textrm{ above } \gamma\} \\
 & -\sharp\{ \textrm{removable }\ i-\textrm{nodes of } \underline{\mu}
\ \textrm{ above } \gamma\},\\
N_i^{b}{(\underline{\lambda},\underline{\mu})}= &   \sharp\{ \textrm{addable }
 i{}-\textrm{nodes of } \underline{\lambda}\ \textrm{ below } \gamma\}\\
&   -\sharp\{ \textrm{removable } i-\textrm{nodes of } \underline{\mu}
\ \textrm{ below } \gamma\},\\
N_{i}{(\underline{\lambda})} =& \sharp\{ \textrm{addable } 
i-\textrm{nodes of } \underline{\lambda}\}\\
& -\sharp\{ \textrm{removable } i-\textrm{nodes of } \underline{\lambda}\},\\
N_{\mathfrak{d}}{(\underline{\lambda})} =& \sharp\{ 0-\textrm{nodes of }
 \underline{\lambda}\}.
\end{align*}

For  $n\in{\mathbb{N}}$, let $\mathcal{F}_n$ be the associated space 
which is defined in the previous paragraph. Let $\mathcal{F}:=
\bigoplus_{n\in{\mathbb{N}}}\mathcal{F}_n $. $\mathcal{F}$ is called
 the Fock space. For each $d$-partition $\underline{\lambda}$, we 
identify the Specht module $S^{\underline{\lambda}}$
 with  $\underline{\lambda}$ so that a basis of $\mathcal{F}_n$ is 
given by the set $\Pi_n^d$.  $\mathcal{F}$ becomes a $\mathcal{U}_q$-
module with the following action:
$$e_{i}\underline{\lambda}=\sum_{\textrm{res}([\underline{\lambda}]
/[\underline{\mu}])=i}{q^{-N_i^{a}{(\underline{\mu},\underline{\lambda})}}
\underline{\mu}},\qquad{f_{i}\underline{\lambda}=\sum_{\textrm{res}
([\underline{\mu}]/[\underline{\lambda}])=i}{q^{N_i^{b}{(\underline{\lambda}
,\underline{\mu})}}\underline{\mu}}},$$
$$k_{h_i}\underline{\lambda}=q^{N_i{(\underline{\lambda})}}\underline{\lambda}
,\qquad{k_{\mathfrak{d}}\underline{\lambda}=q^{-N_{\mathfrak{d}}
{(\underline{\lambda})}}\underline{\lambda}},$$
where $i=0,...,e-1$.  This action was discovered by Hayashi in \cite{Haya}.

Let $\mathcal{M}$ be the $\mathcal{U}_q$-submodule of  $\mathcal{F}$
 generated by the empty $d$-partition. It is isomorphic to an 
integrable  highest weight module (see \cite{DJK}). Now, this result
 allows us to apply the canonical basis theory and the crystal 
graph theory to $\mathcal{M}$.

 In particular, the crystal graph gives a way for labeling the 
Kashiwara-Lusztig's canonical basis of $\mathcal{M}$ (see \cite{Arilivre}
 for more details).   

Based on Misra and Miwa's result,  Ariki and Mathas observed that this 
graph is given by:
\begin{itemize}
\item vertices: the  Kleshchev $d$-partitions,
\item edges: $\displaystyle{{\underline{\lambda}\overset{i}
{\rightarrow}{{\underline{\mu}}}}}$ if and only if $[\underline{\mu}]
/[\underline{\lambda}]$ is a good  $i$-node.
\end{itemize}

Thus, the canonical basis $\mathfrak{B}$ of $\mathcal{M}$ is labeled 
by  the Kleshchev $d$-partitions:

$$\mathfrak{B}=\{ G(\underline{\lambda})\ |\ \underline{\lambda}
\in{\Lambda^0_{\{e;v_0,...,v_{d-1}\}}}\}.$$
 
This set is a basis of  the $\mathcal{U}_{\mathcal{A}}$-module 
$\mathcal{M}_{\mathcal{A}}$ generated by the empty $d$-partition 
and for any specialization  of $q$ into an invertible element ${u}$  
of  a field $R$, we obtain a basis of the specialized module  
$\mathcal{M}_{R,u}$   by specializing the set  $\mathfrak{B}$.

 Now, we have the following theorem of Ariki which shows that 
the problem of computing the decomposition numbers of  $\mathcal{H}_{R,n}$
 can be translated to that of computing the canonical basis of $\mathcal{M}$.
 This theorem was first conjectured by Lascoux, Leclerc and Thibon 
(\cite{LLT}) in the case of  Hecke algebras of type $A_{n-1}$.

\begin{thm}[Ariki \cite{Ari4}]\label{Ar} We have $\Phi_n^d=
\Lambda^0_{\{e;v_0,...,v_{d-1}\}}$. Moreover, assume that $R$ 
is a field of characteristic $0$. Then, for each $\underline{\lambda}
\in{\Phi_n^d}$, there exist polynomials $d_{\underline{\mu},
\underline{\lambda}}(q)\in{\mathbb{Z}[q]}$ and   a unique element 
 $G(\underline{\lambda})$ of the canonical basis such that:
$$G(\underline{\lambda})=\sum_{\underline{\mu}\in{\Pi_n^d}}
{d_{\underline{\mu},\underline{\lambda}}(q)\underline{\mu}}
\qquad{and}\qquad{G(\underline{\lambda})=\underline{\lambda}\ (mod\ q)}.$$
Finally, for all  $\underline{\mu}\in{\Pi_n^d}$ and $\underline{\lambda}
\in{\Phi_n^d}$, we have  $d_{\underline{\mu},\underline{\lambda}}(1)
=d_{\underline{\mu},\underline{\lambda}}$.
\end{thm}

Then, if $R$ is a field of characteristic $0$ and  if we identify 
the $d$-partitions $\underline{\lambda}$ with the modules 
$S^{\underline{\lambda}}$, we see that the canonical basis 
elements specialized at $q=1$ corresponds to the  indecomposable
 projective $\mathcal{H}_{R,n}$-modules.

As noted in the introduction, the problem of this parametrization 
of the simple $\mathcal{H}_{R,n}$-modules  is that we only know 
 a recursive description of the Kleshchev $d$-partitions.  We  
now deal with another parametrization of this set found by Foda 
et  al. which uses almost the same objects as Ariki and Mathas.
\subsection{Parametrization of the simple modules by 
Foda et  al.}\label{FLOTW}

The principal idea of  Foda et al. (\cite{FLOTW}) is to use  
another structure of   $\mathcal{U}_q$-module over $\mathcal{F}$ 
by  choosing another order on the  set of the nodes of the 
$d$-partitions.

 Here, we say that  $\gamma=(a,b,c)$ is above  $\gamma'=(a',b',c')$
 if:
$$b-a+v_c<b'-a'+v_{c'}\ \textrm{or } \textrm{if}\ b-a+v_c=b'-a'+v_{c'}
\textrm{ and }c>c'.$$
This order will be called the FLOTW order.

This order allows us to define functions $\overline{N}_i^{a}
{(\underline{\lambda},\underline{\mu})}$ and  $\overline{N}_i^{b}
{(\underline{\lambda},\underline{\mu})}$ given by the same way as
 ${N}_i^{a}{(\underline{\lambda},\underline{\mu})}$ et    
 ${N}_i^{b}{(\underline{\lambda},\underline{\mu})}$ for the AM order.

Now, we have the following result:

\begin{thm}[Jimbo, Misra, Miwa, Okado \cite{JMMO}] $\mathcal{F}$ 
is a  $\mathcal{U}_q$-module with action:
$$e_{i}\underline{\lambda}=\sum_{\operatorname{res}
([\underline{\lambda}]/[\underline{\mu}])=i}{q^{-\overline{N}_i^{a}
{(\underline{\mu},\underline{\lambda})}}\underline{\mu}},
\qquad{f_{i}\underline{\lambda}=\sum_{\operatorname{res}
([\underline{\mu}]/[\underline{\lambda}])=i}
{q^{\overline{N}_i^{b}{(\underline{\lambda},\underline{\mu})}}
\underline{\mu}}},$$
$$k_{h_i}\lambda=v^{{N}_i{(\underline{\lambda})}}\underline{\lambda}
,\qquad{k_{\mathfrak{d}}\underline{\lambda}=q^{-N_{\mathfrak{d}}
{(\underline{\lambda})}}\underline{\lambda}},$$
 where $0\leq{i}\leq{n-1}$. This action will be called the JMMO action.
\end{thm}

We denote by $\overline{\mathcal{M}}$ the $\mathcal{U}_q$-module 
generated by the empty $d$-partition with the above action. 
This is a highest weight module 
and  the  $d$-partitions of the crystal graph are obtained  
recursively by adding good nodes to  $d$-partitions of the 
crystal graph.

Foda et al. showed that the analogue of the notion of Kleshchev 
$d$-partitions for this action is as follows: 

\begin{df}
 We say that  $\underline\lambda={(\lambda^{(0)} ,...,\lambda^{(d-1)})}$
is a FLOTW  $d$-partition associated to  the set  
${\{e;{v_0},...,v_{d-1}\}}$ if and only if:
\begin{enumerate}
\item for all $0\leq{j}\leq{d-2}$ and $i=1,2,...$, we have:
\begin{align*}
&\lambda_i^{(j)}\geq{\lambda^{(j+1)}_{i+v_{j+1}-v_j}},\\
&\lambda^{(d-1)}_i\geq{\lambda^{(0)}_{i+e+v_0-v_{d-1}}};
\end{align*}
\item  for all  $k>0$, among the residues appearing at the 
right ends of the length $k$ rows of   $\underline\lambda$,
 at least one element of  $\{0,1,...,e-1\}$ does not occur.
\end{enumerate}
We denote by $\Lambda^1_{\{e;{v_0},...,{v_{d-1}}\}}$ the set 
of FLOTW $d$-partitions associated to the set ${\{e;{v_0},...,
{v_{d-1}}\}}$. If there is no ambiguity concerning  
$\{e;{v_0},...,{v_{d-1}}\}$, we denote it by $\Lambda^1$.
\end{df}

Hence, the crystal graph of   $\overline{\mathcal{M}}$ is given by:
\begin{itemize}
\item   vertices: the FLOTW $d$-partitions,
\item  edges: $\displaystyle{{\underline{\lambda}\overset{i}
{\rightarrow}{{\underline{\mu}}}}}$ if and only if  
$[\underline{\mu}]/[\underline{\lambda}]$ is  good  $i$-node 
with respect  to the FLOTW order.
\end{itemize}

So, the canonical basis elements are labeled by the  FLOTW 
$d$-partitions and if we specialize these elements to $q=1$, 
we obtain the same elements as  in Theorem \ref{Ar}. Hence, 
there is a bijection between the set of Kleshchev $d$-partitions
 and the set of FLOTW $d$-partitions.

\section{Schur elements and $a$-functions}

The aim of this section is to introduce Schur elements and 
$a$-functions associated to simple modules of Ariki-Koike algebras. 
It is convenient 
 to express them in terms of symbols. Most of the results explained 
here are given in \cite{BrouKim}. First, we give some notations and 
some definitions:

\subsection{Symbols}\label{symbole}

Let $n\in{\mathbb{N}}$ and $d\in{\mathbb{N}_{>0}}$. The notion of 
symbols are usually associated to $d$-partitions. Here, we generalize
 it to $d$-compositions. This generalization is  useful in our 
argument in  the next section. A $d$-composition $\underline{\lambda}$
 of rank $n$ is a $d$-tuple $(\lambda^{(0)},....,\lambda^{(d-1)})$ where :
\begin{itemize}
\item for all $i=0,...,d-1$, we have $\lambda^{(i)}=(\lambda^{(i)}_1
,...,\lambda^{(i)}_{h^{(i)}})$ for  $h^{(i)}\in{\mathbb{N}}$ and 
$\lambda^{(i)}_j\in{\mathbb{N}_{>0}}$ ($j=1,...,h^{(i)}$), $h^{(i)}$ 
is called the height  of $\lambda^{(i)}$,
\item $\displaystyle{\sum_{i=0}^{d-1}\sum_{j=1}^{h^{(i)}}
\lambda^{(i)}_j=n}$.
\end{itemize}

 Let $\underline{\lambda}=(\lambda^{(0)},...,\lambda^{(d-1)})$ 
be a $d$-composition and let $h^{(i)}$ be the heights of the
 compositions $\lambda^{(i)}$. Then the height of 
$\underline{\lambda}$ is the following positive integer:
$$h_{\underline{\lambda}}=\textrm{max}{\{h^{(0)},...,h^{(d-1)}\}}$$

Let $k$ be a positive integer. The ordinary symbol associated to 
$\underline{\lambda}$ and $k$ is the following set:
$$\textbf{B}:=(B^{(0)},...,B^{(d-1)}),$$
where $B^{(i)}$, for  $i=0,...,d-1$, is given by:
$$B^{(i)}:=(B^{(i)}_1,...,B^{(i)}_{h_{\underline{\lambda}}+k}),$$
in which:
$$B^{(i)}_{j}:=\lambda^{(i)}_j-j+h_{\underline{\lambda}}+k\ \ 
(1\leq j \leq h_{\underline{\lambda}}+k),$$
and $\lambda^{(i)}_j:=0$ if $j>h^{(i)}$. Note that the entries 
of each $B^{(i)}$ ($i=0,...,d-1$) are given in descending order, 
contrary to the usual convention.

We denote by $h_{\textbf{B}}:=h_{\underline{\lambda}}+k$ the 
height of the symbol $\textbf{B}$. Now, we fix a sequence of 
rational positive numbers:
$$m=(m^{(0)},...,m^{(d-1)}).$$
Let $k$ be a positive integer. We define a set $\textbf{B}[m]'$ 
associated to $m$ and $k$ by adding on each part $B^{(i)}_j$ of the 
ordinary symbol the number $m^{(i)}$.
So, we have:
$$\textbf{B}[m]'=(B'^{(0)},...,B'^{(d-1)}),$$
where $B'^{(i)}$, for $i=0,...,d-1$, is given by:
$$B'^{(i)}:=(B'^{(i)}_1,...,B'^{(i)}_{h_{\underline{\lambda}}+k}),$$
in which:
$$B'^{(i)}_{j}:=\lambda^{(i)}_j-j+h_{\underline{\lambda}}+k+m^{(i)}
\  \  (1\leq j \leq h_{\underline{\lambda}}+k).$$ 
Note that this set is defined in \cite{BrouKim} when the numbers 
$m^{(i)}$ are integers using the notion of $m$-translated symbols.
\begin{exa}  
Let $d=3$ and $m=(1,\displaystyle{\frac{1}{2}},2)$. We consider 
the following $3$-partition:
$$\underline{\lambda}=((4,2),(0),(5,2,1)).$$
Let $k=0$. The ordinary symbol associated to $\underline{\lambda}$ 
is given by:
$$\textbf{B}=\left\lbrace
 \begin{array}{ccccc}
   B^{(0)} & = & 6  & 3 & 0 \\
    B^{(1)} & = & 2  & 1 & 0 \\
    B^{(2)} & = & 7  & 3 & 1 \\
\end{array}
\right.$$
The set  $\textbf{B}[m]'$ associated to $\underline{\lambda}$ is given by:
$$\textbf{B}[m]'=\left\lbrace
 \begin{array}{ccccc}
   B'^{(0)} & = & 7  & 4 & 1 \\
    B'^{(1)} & = & {5}/{2}  &{3}/{2} &{1}/{2} \\
    B'^{(2)} & = & 9  & 5 & 3 \\
\end{array}
\right.$$
\end{exa}

\subsection{Schur elements}
Let $\mathcal{H}_{R,n}$ be a semi-simple Ariki-Koike algebra over 
a field $R$ with parameters $\{v;u_0,...,u_{d-1}\}$.
 $\mathcal{H}_{R,n}$ is a symmetric algebra and thus, it has a 
symmetrizing normalized trace $\tau$ (see \cite{BrouKim} and 
\cite[chapter 7]{geck1} for the representation theory of symmetric 
algebras). The Schur elements are defined to be the non zero 
elements $s_{\underline{\lambda}}(v,u_0,...,u_{d-1})\in{R}$ such that:
$$\tau=\sum_{\underline{\lambda}\in{\Pi_{n}^d}}  \frac{1}
{s_{\underline{\lambda}}(v,u_0,...,u_{d-1})}\chi_{\underline{\lambda}},$$
where, for $\underline{\lambda}\in{\Pi_{n}^d}$, 
$\chi_{\underline{\lambda}}$ is the  irreducible character 
associated to the simple  $\mathcal{H}_{R,n}$-module $S^{\underline{\lambda}}$.

Now, we give the expressions of these Schur elements. They were 
 conjectured in \cite{Ma}, calculated in \cite{GIM}  and the 
following expression is given in  \cite{BrouKim}.

First, let's fix some notations:
 let $\textbf{B}=(B^{(0)},...,B^{(d-1)})$ be an ordinary  symbol 
of height $h$ associated to  a  $d$-partition $\underline{\lambda}=(
\lambda^{(0)},...,\lambda^{(d-1)})$ of rank $n$. Following 
\cite{BrouKim},  we denote:

\begin{align*}
\delta_{\textbf{B}}(v,u_0,...,u_{d-1})&=\prod_{{0\leq{i}\leq{j}<d}
\atop{(\alpha,\beta)\in{B^{(i)}\times{B^{(j)}}}\atop{\alpha>\beta\
 \textrm{if}\ i=j}}}(v^{\alpha}u_i-v^{\beta}u_j),\\  
\theta_{\textbf{B}}(v,u_0,...,u_{d-1})&=\prod_{{0\leq{i,j}<d}
\atop{\alpha\in{B^{(i)}}\atop{1\leq{k}\leq{\alpha}}}}(v^{k}u_i-u_j),\\
\sigma_{\textbf{B}}&=\left( \begin{array}{c} d \\ 2 \end{array} \right)
  \left( \begin{array}{c} h \\ 2 \end{array} \right)+ n(d-1),\\
 \tau_{\textbf{B}}&= \left( \begin{array}{c} d(h-1)+1 \\ 2 
\end{array} \right)  + \left( \begin{array}{c} d(h-2)+1 \\ 2 
\end{array} \right)+...+   \left( \begin{array}{c} 2 \\ 2 \end{array} 
\right),\\ |\textbf{B}|&=\sum_{{0\leq{i}<d}\atop{\alpha\in{B^{(i)}}}}\alpha,
 \\
\nu_{\textbf{B}}(v,u_0,...,u_{d-1})&=\prod_{0\leq{i}<j<d}
(u_i-u_j)^{h}\theta_{\textbf{B}}(v,u_0,...,u_{d-1}). \end{align*}
Now, we have the following proposition:
\begin{prp}[Geck, Iancu, Malle, see \cite{BrouKim}]\label{afonction} 
Let $S^{\underline{\lambda}}$ be the simple $\mathcal{H}_{R,n}$-module
 associated to the  $d$-partition $\underline{\lambda}$  and let  
$\textbf{B}$ be an ordinary symbol   associated to   $\underline{\lambda}$,
 then the Schur element of  $S^{\underline{\lambda}}$ is given by:
$$s_{\underline{\lambda}}(v,u_0,...,u_{d-1})=\left( (v-1)
\prod_{0\leq{i}<d}u_i\right)^{-n}(-1)^{\sigma_{\textbf{B}}}
v^{\tau_{\textbf{B}}-|\textbf{B}|+n}\frac{\nu_{\textbf{B}}
(v,u_0,...,u_{d-1})}{\delta_{\textbf{B}}(v,u_0,...,u_{d-1})}.$$
\end{prp}

\subsection{$a$-function for  Ariki-Koike algebras}\label{princip}

In this paragraph, we deal with the   $a$-value  of  the  simple  
modules in the case of  Ariki-Koike  algebras. 
 Let $d$ and $e$ be two positive integers and define complex numbers by:
$$\eta_d=\textrm{exp}({\frac{2i\pi}{d}})\qquad \textrm{and}
\qquad{\eta_e=\textrm{exp}({\frac{2i\pi}{e}})}.$$

 We fix a sequence of rational numbers   $m=(m^{(0)},...,m^{(d-1)})$ 
such that  we have  $dm^{(j)}\in{\mathbb{N}}$ for all $j=0,...,d-1$.  
Let  $y$ be an indeterminate and let $L:=\mathbb{Q}[\eta_d](y)$.  
We  consider the Ariki-Koike algebra $\mathcal{H}_{L,n}$ defined 
over the field $L$ with the following choice of parameters:
\begin{eqnarray*}
&&u_j=y^{dm^{(j)}}\eta_d^{j} \textrm{ for}\ j=0,...,d-1,\\
&&v=y^d. 
\end{eqnarray*}
 By  Theorem \ref{semisimple}, $\mathcal{H}_{L,n}$ is a split 
semi-simple Ariki-Koike algebra.Thus,  the simple modules are 
 labeled by the $d$-partitions of rank $n$.

We can now define the $a$-value of  the simple $\mathcal{H}_{L,n}$-modules. 
Let $\underline{\lambda}$ be a $d$-partition of rank $n$ and let 
  $S^{{\underline{\lambda}}}$ be the simple module associated to 
this $d$-partition. Let $s_{\underline{\lambda}}$ be the Schur 
element of $S^{{\underline{\lambda}}}$. This is a Laurent polynomial in $y$. 
Let $\textrm{val}_y{(s_{\underline{\lambda}})}$ be the valuation 
of  $s_{\underline{\lambda}}$ in $y$ that is to say:
 $$\textrm{val}_y{(s_{\underline{\lambda}})}:=
{}-\textrm{min}{\{l\in{\mathbb{Z}}\ |\ y^ls_{\underline{\lambda}}
\in{\mathbb{Z}[\eta_d]}[y]\}} $$
Then, the $a$-value of  $S^{{\underline{\lambda}}}$ is defined by:
$$a(S^{\underline{\lambda}}):=-\frac{\textrm{val}_y
{(s_{\underline{\lambda}})}}{d}.$$
For all $\underline{\lambda}\in{\Pi_n^d}$,  $a(S^{\underline{\lambda}})$
 is a rational number which depends on  $(m^{(0)},...,m^{(d-1)})$.
 To simplify, we will denote  $a(\underline{\lambda}):= 
a(S^{\underline{\lambda}})$.

\begin{prp}\label{af} Let  $\underline{\lambda}$ be a $d$-partition 
of rank $n$  and let  $S^{\underline{\lambda}}$ be the simple 
$\mathcal{H}_{L,n}$-module.

Let  $k$ be a positive integer and $\textbf{B}[m]'$ be the set 
associated to $m$ and $k$  as in paragraph \ref{symbole}. Then, 
if $h$ is the height of  $\textbf{B}[m]'$, we have:
$$a(\underline{\lambda})=f(n,h,m)+\sum_{{0\leq{i}\leq{j}<d}
\atop{{(a,b)\in{B'^{(i)}\times{B'^{(j)}}}}\atop{a>b\ 
\textrm{if}\ i=j}}}{\min{\{a,b\}}}   -\sum_{{0\leq{i,j}<d}
\atop{{a\in{B'^{(i)}}}\atop{1\leq{k}\leq{a}}}}{\min{\{k,m^{(j)}\}}},$$
where:
\begin{eqnarray*}
&f(n,h,m)&\displaystyle{=n\sum_{j=0}^{d-1}{m^{(j)}}-
\tau_{\textbf{B}}+|\textbf{B}|-n-h\sum_{0\leq{i}<j<d}
\min{\{m^{(i)},m^{(j)}\}}+}\\
 &&\displaystyle{+\sum_{{0\leq{i,j}<d}\atop{\alpha\in{B^{(i)}}
\atop{1\leq{k}\leq{m^{(i)}}}}}{\min{\{k,m^{(j)}\}}}}.
\end{eqnarray*}
\end{prp}
\begin{proof} 
We use Proposition \ref{afonction}:
$$a(\underline{\lambda})=n\sum_{j=0}^{d-1}{m^{(j)}}-
\tau_{\textbf{B}}+|\textbf{B}|-n-h\sum_{0\leq{i}<j<d}
\min{\{m^{(i)},m^{(j)}\}}+$$
$$+\sum_{{0\leq{i}\leq{j}<d}\atop{{(\alpha,\beta)\in{B^{(i)}
\times{B^{(j)}}}}\atop{\alpha>\beta\ \textrm{if}\ i=j}}}
{\min{\{\alpha+m^{(i)},\beta+m^{(j)}\}}}-\sum_{{0\leq{i,j}<d}
\atop{\alpha\in{B^{(i)}}\atop{1\leq{k}\leq{\alpha}}}}
{\min{\{k+m^{(i)},m^{(j)}\}}}.$$
We have:
$$\sum_{{0\leq{i,j}<d}\atop{\alpha\in{B^{(i)}}\atop{1\leq{k}
\leq{\alpha}}}}{\min{\{k+m^{(i)},m^{(j)}\}}}=\sum_{{0\leq{i,j}<d}
\atop{\alpha\in{B^{(i)}}\atop{1\leq{k}\leq{\alpha+m^{(i)}}}}}
{\min{\{k,m^{(j)}\}}}-\sum_{{0\leq{i,j}<d}\atop{\alpha\in{B^{(i)}}
\atop{1\leq{k}\leq{m^{(i)}}}}}{\min{\{k,m^{(j)}\}}}.$$
Using the set $\textbf{B}[m]'$ of paragraph \ref{symbole}, we conclude that:
$$a(\underline{\lambda})=f(n,h,m)+\sum_{{0\leq{i}\leq{j}<d}
\atop{{(a,b)\in{B'^{(i)}\times{B'^{(j)}}}}\atop{a>b\ \textrm{if}\ i=j}}}
{\min{\{a,b\}}}   -\sum_{{0\leq{i,j}<d}\atop{{a\in{B'^{(i)}}}
\atop{1\leq{k}\leq{a}}}}{\min{\{k,m^{(j)}\}}}.$$\end{proof}
\begin{re}
The formula of the $a$-function does not depend on the choice of
 the ordinary symbol.
\end{re}

The next section gives an interpretation of the parametrization of 
the simple modules   by Foda et al. in terms of this $a$-function when
 the Ariki-Koike algebra is not semi-simple.

\section{Foda et al. $d$-partitions and $a$-functions}

First, we introduce some definitions:

\subsection{Notations and  hypothesis}\label{not}

Here, we keep the notations of the previous section. As in the 
second section, we assume that we have nonnegative integers :
$$0\leq{v_0}\leq v_1 \leq ... \leq v_{d-1}<e.$$
Then,  for $j=0,...,d-1$, we define rational numbers:
$$m^{(j)}=v_j-\frac{je}{d}+se,$$
where $s$ is a positive  integer such that $m^{(j)}\geq{0}$ 
 for $j=0,...,d-1$. We have  $dm^{(j)}\in{\mathbb{N}}$ for 
all $j=0,...,d-1$.

Let   $\mathcal{H}_{L,n}$ be the Ariki-Koike algebra over
 $L=\mathbb{Q}[\eta_d](y)$  with the following parameters:
\begin{eqnarray*}
&&u_j=y^{dm^{(j)}}\eta_d^{j}\ \textrm{for}\ j=0,...,d-1,\\
&&v=y^d. 
\end{eqnarray*}

If we specialize the parameter $y$ into ${\eta_{de}:=\textrm{exp}
(\frac{2i\pi}{de})}$, we obtain an Ariki-Koike algebra  $\mathcal{H}_{R,n}$
 over $R=\mathbb{Q}(\eta_{de})$ with the following parameters:
\begin{eqnarray*}
&&u_j=\eta_{de}^{dm^{(j)}}\eta_d^{j}=\eta_e^{v_j}\ 
\textrm{for}\ j=0,...,d-1,\\
&&v=\eta_e. 
\end{eqnarray*}
This algebra is split but  not  semi-simple and we will discuss 
about its representation theory.

 Here, we keep the notations used at the previous section and we 
add the following ones.
First, we extend the combinatorial definitions of diagrams and  
residues introduced in  paragraph \ref{combinatoire} to the $d$-
compositions of rank $n$ in an obvious way.  
Let $\underline{\lambda}$ be a $d$-composition of rank $n$. 
The nodes at the right ends of the parts of $\underline{\lambda}$ 
will be  called the nodes of the border of $\underline{\lambda}$.
\begin{exa} 
Assume that $d=2$, $v_0=0$, $v_1=2$ and that  $e=4$. Let 
$\underline{\lambda}=(4.2.3,3.5)$, then the diagram of  
$\underline{\lambda}$ is the following one:

$$  \underline{\lambda}  =\left( \ \begin{tabular}{|c|c|c|c|}
                         \hline
                           0  & 1  & 2  & \textbf{3}  \\
                          \hline
                            3  & \textbf{0}\\
                          \cline{1-3}
                         2  & 3 & \textbf{0}\\
  \cline{1-3}
                             \end{tabular}\  ,\
                       \begin{tabular}{|c|c|c|c|c|}
                         \cline{1-3}
                           2  & 3  & \textbf{0}   \\
                          \hline
                          1 & 2 & 3 & 0 &     \textbf{ 1 }    \\
                          \hline
                             \end{tabular}\
                                                 \right)$$
\\
The nodes in bold-faced type are the nodes of the border 
of $\underline{\lambda}$.
\end{exa}

We will also use  the following notations: Let  $\underline{\lambda}$ 
and $\underline{\mu}$ be   $d$-compositions  of rank  
$n$ and $n+1$ and let  $k\in{\{0,...,e-1\}}$ be a residue.
\begin{itemize}
\item We write $\underline{\lambda}{\overset{k}{\underset{(j,p)}
{\mapsto}}}{\underline{\mu}}$ if  $[\underline{\mu}]/ 
[\underline{\lambda}]$ is  a  $k$-node  on the  part $\mu^{(p)}_j$.
\item We write  $\underline{\lambda}{\overset{k}{\mapsto}
{\underline{\mu}}}$ if  $[\underline{\mu}] / [\underline{\lambda}]$
 is  a  $k$-node.
\end{itemize}

The next paragraph gives combinatorial properties concerning FLOTW 
$d$-partitions.

\subsection{Combinatorial properties}
First, we give some lemmas:
\begin{lmm}\label{premierlemme} Let $\underline{\lambda}$ be a  
FLOTW $d$-partition  associated to  $\{e;v_{0},...,v_{d-1}\}$. 
Let  $\xi$ be a  removable node on a  part $\lambda_{j_1}^{(i_1)}$, 
let $\lambda_{j_2}^{(i_2)}$ be a part of  $\underline{\lambda}$ 
and  assume that:
$$\lambda_{j_2}^{(i_2)}-j_2+v_{i_2}\equiv{\lambda_{j_1}^{(i_1)}
{}-j_1+v_{i_1}-1\ (\textrm{mod }e)}.$$
  Then, we have:

$$\lambda^{(i_2)}_{j_2}\geq{\lambda_{j_1}^{(i_1)}}\iff{\lambda_
{j_2}^{(i_2)}-j_2+m^{(i_2)}+1\geq{\lambda_{j_1}^{(i_1)}-j_1+m^{(i_1)}}}.$$
\end{lmm}
\begin{proof} 
Using the hypothesis, there exists some  $t\in{\mathbb{Z}}$ such that:
$$\lambda_{j_2}^{(i_2)}-j_2+v_{i_2}=\lambda_{j_1}^{(i_1)}-j_1+v_{i_1}-1+te.$$
a) We assume that:
$$\lambda_{j_2}^{(i_2)}-j_2+m^{(i_2)}+1\geq{\lambda_{j_1}^{(i_1)}-
j_1+m^{(i_1)}}.$$
We want to show:
$$\lambda^{(i_2)}_{j_2}\geq{\lambda_{j_1}^{(i_1)}}.$$
Assume to the contrary that:
$$\lambda^{(i_2)}_{j_2}<{\lambda_{j_1}^{(i_1)}}.$$
We have:
$$\lambda_{j_2}^{(i_2)}-j_2+v_{i_2}-i_2\frac{e}{d}+1
\geq{\lambda_{j_1}^{(i_1)}-j_1+v_{i_1}-i_1\frac{e}{d}},$$
then:
$$-1+te-i_2\frac{e}{d}+1\geq{-i_1\frac{e}{d}}.$$
So, we have:
$$te\geq{(i_2-i_1)\frac{e}{d}}.$$
We now distinguish  two cases.\\
\underline{If $i_1\geq{i_2}$}:\\
Then, we have $t\geq{0}$ because $i_2-i_1\geq 1-d$, thus:
$$\lambda_{j_2}^{(i_2)}-j_2+v_{i_2}\geq{\lambda_{j_1}^{(i_1)}
{}-j_1+v_{i_1}-1},$$
so:
$$j_1-j_2\geq{v_{i_1}-v_{i_2}}.$$
Now, we use the characterization of the FLOTW $d$-partitions:
$$\lambda_{j_2}^{(i_2)}\geq{\lambda^{(i_1)}_{j_2+v_{i_1}-v_{i_2}}}.$$
We obtain:
$$\lambda_{j_2}^{(i_2)}\geq{\lambda_{j_1}^{(i_1)}},$$
a contradiction.\\
\\
\underline{If $i_1<i_2$}:\\
Then, we have  $t>0$, thus:
$$\lambda_{j_2}^{(i_2)}-j_2+v_{i_2}\geq{\lambda_{j_1}^{(i_1)}
{}-j_1+v_{i_1}-1+e},$$
so:
$$j_1-j_2\geq{v_{i_1}-v_{i_2}+e}.$$
 We use the characterization of the FLOTW $d$-partitions:
$$\lambda_{j_2}^{(i_2)}\geq{\lambda^{(i_1)}_{j_2+v_{i_1}-v_{i_2}+e}}.$$
We obtain:
$$\lambda_{j_2}^{(i_2)}\geq{\lambda_{j_1}^{(i_1)}},$$
a contradiction again.\\
\\
b) We assume that:
$$\lambda_{j_2}^{(i_2)}-j_2+m^{(i_2)}+1<{\lambda_{j_1}^{(i_1)}
{}-j_1+m^{(i_1)}}.$$
We want to show:
$$\lambda^{(i_2)}_{j_2}<{\lambda_{j_1}^{(i_1)}}.$$
Assume to the contrary that:
$$\lambda^{(i_2)}_{j_2}\geq{{\lambda_{j_1}^{(i_1)}}}.$$
As above, we have:
$$te<(i_2-i_1)\frac{e}{d}.$$
Again, there are two cases to consider.\\
\underline{If $i_1>i_2$}:\\
Then, we have $t<{0}$, and so:
$$\lambda^{(i_2)}_{j_2}-j_2+v_{i_2}\leq{\lambda^{(i_1)}_
{j_1}-j_1+v_{i_1}-1-e},$$
thus:
$$j_2-(j_1+1)\geq{v_{i_2}-v_{i_1}+e}.$$
We use the characterization of the FLOTW $d$-partitions:
$$\lambda_{j_1+1}^{(i_1)}\geq{\lambda^{(i_2)}_{j_1+1+e+
v_{i_2}-v_{i_1}}},$$
so:
$$\lambda_{j_1+1}^{(i_1)}\geq{\lambda^{(i_2)}_{j_2}}.$$
Now,  $\xi$ is a removable node of $\underline{\lambda}$. 
Thus, we obtain:
$$\lambda_{j_1}^{(i_1)}>\lambda_{j_1+1}^{(i_1)},$$
a contradiction.\\
\underline{If $i_1\leq{i_2}$}:\\
Then, we have  $t\leq{0}$ because $i_1 -i_2 \geq 1-d$, thus:
$$\lambda^{(i_2)}_{j_2}-j_2+v_{i_2}\leq{\lambda^{(i_1)}_{j_1}
{}-j_1+v_{i_1}-1},$$
so:
$$j_2-(j_1+1)\geq{v_{i_2}-v_{i_1}}.$$
We use the characterization of the FLOTW $d$-partitions:
$$\lambda_{j_1+1}^{(i_1)}\geq{\lambda^{(i_2)}_{j_1+1+
v_{i_2}-v_{i_1}}},$$
so:
$$\lambda_{j_1+1}^{(i_1)}\geq{\lambda^{(i_2)}_{j_2}}.$$
Now,  $\xi$ is a removable node of $\underline{\lambda}$,
 so:
$$\lambda_{j_1}^{(i_1)}>\lambda_{j_1+1}^{(i_1)}.$$
Again, this is a contradiction.
\end{proof}

\begin{lmm}\label{deuxiemelemme} Let $\underline{\lambda}$ be 
a FLOTW $d$-partition associated to  $\{e;v_{0},...,v_{d-1}\}$. We define:
 $$l_{\textrm{max}}:=\operatorname{max}\{\lambda^{(0)}_1 ,
\lambda^{(1)}_1,..,, \lambda^{(l-1)}_1\}$$
Then, there exists a removable $k$-node $\xi_1$, for some $k$, 
 on a part $\lambda_{j_1}^{(i_1)}=l_{max}$   which satisfies 
the following property:
 if $\xi_2$ is a  $k-1$-node on the border of a part  
$\lambda_{j_2}^{(i_2)}$, then:
$$\lambda^{(i_1)}_{j_1}>\lambda^{(i_2)}_{j_2}.$$
\end{lmm}
\begin{proof} 
Let  $\lambda^{(l_1)}$, ..., $\lambda^{(l_r)}$ be the partitions 
of $\underline{\lambda}$ such that  $\lambda_1^{(l_1)}=...=
\lambda_1^{(l_r)}=l_{\textrm{max}}$ are the parts of  maximal length. 
 Let  $k_1,...,k_r$ be the residues of the removable nodes  
$\xi_1$,...,$\xi_r$ on the parts of length  $l_{\textrm{max}}$. 
We want  to show that there exists  $1\leq{i}\leq{r}$ such that 
there is no node with  residue $k_i-1$ on the border of a part with length $l_{\textrm{max}}  $. 

Assume that, for each $1\leq{i}\leq{r}$, there exists a node 
on the border of a part of length   $l_{\textrm{max}}$ with 
residue  $k_i-1$. 

 Then, there exists a partition  $\lambda^{(l_{s_1})}$, for 
some $1\leq s_1 \leq r$,  with a $k_1-1$-node on the border 
of a part  of length $l_{\textrm{max}}$:\\
\\
\\
$$\lambda^{(l_{s_1})}= \ \begin{array}{|ccc|c|}
                        \hline
                        \ \ \      & ...  &\ \ \   & ... \\  
                          \ \ \  & ...  &\ \ \   &  ...  \\
                      \hline
                         \ \ \  &\ \ \     &\ \ \   & k_1-1 \\
                         \hline
              \ \ \  & ...  &\ \ \   & ... \\
                          \ \ \  & ...  &\ \ \   & ...  \\
                      \hline
                         \ \ \  &\ \ \     &\ \ \   & k_{s_1} \\ 
                        \hline
               \ \ \  & ...     &\ \ \  \\
               \ \ \  & ...    &\ \ \  \\
                           \end{array}\  $$
\begin{pspicture}(20,0.4)
\psline{<->}(4.9,0.4)(8.5,0.4)
\rput(6.7,0){$l_{\textrm{max}}$}
\end{pspicture}
We have  $s_1\neq{1}$, otherwise   the nodes on the border of the 
parts with length    $l_{\textrm{max}}$ on  $\lambda^{(l_{s_1})}$ 
would describe the set  $\{0,...,e-1\}$. This violates the  second
 condition to be a  FLOTW $d$-partition.

  We use the same idea for the residue  $k_{s_1}$, there exists 
$\lambda^{(l_{s_2})}$, for some $1\leq s_2\leq r$, as below:

$$\lambda^{(l_{s_2})}= \ \begin{array}{|ccc|c|}
                        \hline
                          \ \ \  & ...  &\ \ \   & ... \\
                          \ \ \  & ...  &\ \ \   &  ...  \\
                      \hline
                         \ \ \  &\ \ \     &\ \ \   & k_{s_1}-1 \\
                         \hline
              \ \ \  & ...  &\ \ \   & ... \\
                          \ \ \  & ...  &\ \ \   & ...  \\ 
                      \hline
                         \ \ \  &\ \ \     &\ \ \   & k_{s_2} \\ 
                        \hline
               \ \ \  & ...     &\ \ \  \\
               \ \ \  & ...    &\ \ \  \\
                           \end{array}\  $$
\begin{pspicture}(20,0.4)
\psline{<->}(4.9,0.4)(8.5,0.4)
\rput(6.7,0){$l_{\textrm{max}}$}
\end{pspicture}\\
\\
We have  $s_2 \neq{s_1}$ (for the same reasons as above) and 
$s_2\neq{1}$, otherwise  the nodes on the border of the parts 
with length    $l_{\textrm{max}}$ on  $\lambda^{(l_{s_2})}$ and 
on $\lambda^{(l_{s_1})}$ would describe all the set  $\{0,...,e-1\}$. 
Continuing in this way, we finally obtain that there exists 
$1\leq s_r\leq r$ such that:

$$s_r\notin{\{1,s_1,s_2,...,s_{r-1}\}}.$$
This is impossible since   $s_r\in{\{1,...,r\}}$ and for all $i\neq j$,
 we have   $s_i\neq s_j$.

So, there exists $0\leq{i}\leq{r}$ such that there is no  $k_i-1$-node
 on the border of the  parts  with maximal length.
\end{proof}

\begin{lmm}\label{troisiemelemme} Let $\underline{\lambda}$ be a FLOTW 
$d$-partition associated to   $\{e;v_{0},...,v_{d-1}\}$. Let  $k$ 
and $\xi_1$ be as in the previous lemma.

Let $\xi_1$, $\xi_2$,..., $\xi_s$ be the removable $k$-nodes  of 
 $\underline{\lambda}$. We assume that they are on parts  
$\displaystyle{\lambda_{j_1}^{(i_1)}\geq{\lambda_{j_2}^{(i_2)}
\geq{...}\geq{\lambda_{j_s}^{(i_s)}}}}$.

Let $\gamma_1$, $\gamma_2$,..., $\gamma_r$ be the  $k-1$-nodes 
on the border of   $\underline{\lambda}$. We assume that they are
  on parts $\displaystyle{\lambda_{p_1}^{(l_1)}\geq{\lambda_{p_2}^{(l_2)}
\geq{...}\geq{\lambda_{p_r}^{(l_r)}}}}$.

We remove the nodes $\xi_u$ such that  $\lambda_{j_u}^{(i_u)}   >
 \lambda_{p_1}^{(l_1)}$ from $\underline{\lambda}$. Let   
$\underline{\lambda}'$ be the  resulting $d$-partition. Then,  
 $\underline{\lambda}'$ is a FLOTW $d$-partition  associated to 
 $\{e;v_{0},...,v_{d-1}\}$ and the  rank of  $\underline{\lambda}'$ 
is   strictly smaller than the rank of  $\underline{\lambda}$.
\end{lmm}
\begin{proof} 
The rank of  $\underline{\lambda}'$ is strictly smaller than the 
rank of  $\underline{\lambda}$ by the previous lemma.
We verify the two conditions of the FLOTW $d$-partitions  for 
 $\underline{\lambda}'$.\\
\\
\underline{First Condition}:\\
 First, we have to verify that if $\lambda_j^{(i)}> 
\lambda_{p_1}^{(l_1)}$,  $\lambda_j^{(i)}=\lambda_{j+v_{i+1}-v_{i}}^{(i+1)}$ 
and if we remove a node from   $\lambda_j^{(i)}$, then, we also remove 
 a node from $\lambda_{j+v_{i+1}-v_{i}}^{(i+1)}$. Observe that:
$$\lambda_{j+v_{i+1}-v_{i}}^{(i+1)}-(j+v_{i+1}-v_{i})+v_{i+1}=
\lambda_{j}^{(i)}-j+v_{i}.$$
Thus,  the residue on the border of   $\lambda_{j+v_{i+1}-v_{i}}^{(i+1)}$
 is  $k$. Moreover, the associated node is a removable one,  otherwise
 we would have:
$$\lambda^{(i+1)}_{j+v_{i+1}-v_{i}+1}>\lambda^{(i)}_{j+1},$$
 contradicting to our assumption that $\underline{\lambda}$ is a 
FLOTW $d$-partition.

Since  $\lambda_j^{(i)}=\lambda_{j+v_{i+1}-v_{i}}^{(i+1)}$ and  
$\lambda_{j+v_{i+1}-v_i}^{(i+1)}> \lambda_{p_1}^{(l_1)}$, the  $k$-node 
on the border of  $\lambda_{j+v_{i+1}-v_{i}}^{(i+1)}$ must be removed.
 So, the first condition of the FLOTW $d$-partitions holds for 
$\underline{\lambda}'$.\\
\\
\underline{Second Condition }:\\
The only problem may arrive when there exists  $t\in{\{1,...,s}\}$ 
such that:
\begin{itemize}
\item if we delete the nodes   $\xi_1$, $\xi_2$,..., $\xi_{t-1}$, 
the resulting $d$-partition $\underline{\lambda}\setminus
{\{\xi_1,...,\xi_{t-1}\}}$ satisfies the second condition.
\item   $\underline{\lambda}\setminus{\{\xi_1,...,\xi_{t}\}}$ 
doesn't satisfy the second condition. 
\end{itemize}
This implies that the set of residues of the nodes on the border 
of the parts of $\underline{\lambda}$ with length  
$\lambda_{j_t}^{(i_t)}-1$   is equal to  the following set:
$$\{0,...,e-1\}\setminus{\{k-1\}}.$$
Note that the second condition is satisfied for all the other 
lengths than  $\lambda_{j_t}^{(i_t)}-1$. 
For example, this problem occurs when    $\underline{\lambda}$ 
is as follows:
\\
\\
\
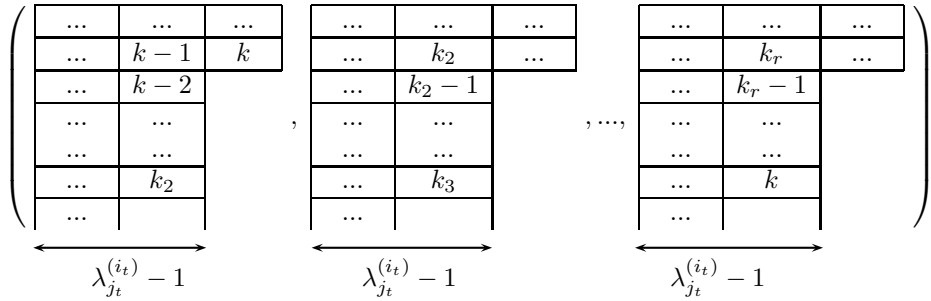
$$\left( \ \begin{array}{|c|c|c|}
                        \hline
                          \ \ ... \ \     & ... & ... \\  
                        \hline
                          \ \   ...\  \   &   k-1   &\ \  k \ \ \\ 
                      \hline
                         \ \  ...\ \      & k-2 \\        
                         \cline{1-2}
              \ \   ...  \ \ & ... \\  
              \ \ ... \ \  & \ \ ... \ \   \\ 
                      \cline{1-2}
                         \ \  ...\ \       & k_2 \\ 
                        \cline{1-2}
               \ \   ...\ \      &  \\
               
                           \end{array}\ 
,\  \begin{array}{|c|c|c|}
                        \hline
                          \ \ ... \ \     & ... & ... \\  
                        \hline
                          \ \   ...\  \   &   k_2   &\ \  ... \ \ \\ 
                      \hline
                         \ \  ...\ \      & k_2-1 \\        
                         \cline{1-2}
              \ \   ...  \ \ & ... \\  
              \ \ ... \ \  & \ \ ... \ \   \\ 
                      \cline{1-2}
                         \ \  ...\ \       & k_3 \\ 
                        \cline{1-2}
               \ \   ...\ \      &  \\
               
                           \end{array}\ 
,...,\  \begin{array}{|c|c|c|}
                        \hline
                          \ \ ... \ \     & ... & ... \\  
                        \hline
                          \ \   ...\  \   &   k_r   &\ \  ... \ \ \\ 
                      \hline
                         \ \  ...\ \      & k_r-1 \\        
                         \cline{1-2}
              \ \   ...  \ \ & ... \\  
              \ \ ... \ \  & \ \ ... \ \   \\ 
                      \cline{1-2}
                         \ \  ...\ \       & k \\ 
                        \cline{1-2}
               \ \   ...\ \      &  \\
               
                           \end{array}\ 
 \right)  $$
\begin{pspicture}(20,0.4)
\psline{<->}(0.4,0.4)(2.7,0.4)
\psline{<->}(4.1,0.4)(6.5,0.4)
\psline{<->}(8.4,0.4)(10.9,0.4)
\rput(1.8,0){$\lambda_{j_t}^{(i_t)}-1$}
\rput(5.4,0){$\lambda_{j_t}^{(i_t)}-1$}
\rput(9.5,0){$\lambda_{j_t}^{(i_t)}-1$}
\end{pspicture}\\
\\
We want to show that among the residues on the border of the parts 
of $\underline{\lambda}'$ with length  $\lambda_{j_t}^{(i_t)}-1$,
 $k$ does not occur.

There exists at least one  $k$-node on the border of a part of 
$\underline{\lambda}$ with length  $\lambda_{j_t}^{(i_t)}-1$. 
Such a $k$-node must be a removable one. If otherwise, we would 
have  a $k-1$-node on the border of a part of $\underline{\lambda}$ 
with length  $\lambda_{j_t}^{(i_t)}-1$ contradicting to our assumption.

 We have $\lambda_{j_t}^{(i_t)}-1>\lambda_{p_1}^{(l_1)}$ because 
there is no  $k{}-1$-node on parts with length  $\lambda_{j_t}^{(i_t)}-1$.
So,  all the $k$-nodes on the border of parts with length 
 $\lambda_{j_1}^{(i_1)}-1$ must be removed.  Thus, we do 
remove all of them.

Then, the set of residues of the nodes on the border of parts
 of $\underline{\lambda}'$ with length  $\lambda_{j_t}^{(i_t)}-1$ 
is equal to:
$$\{0,...,e-1\}\setminus{\{k\}}.$$
Keeping the above example, $\underline{\lambda}'$ is given by:
\\
\\
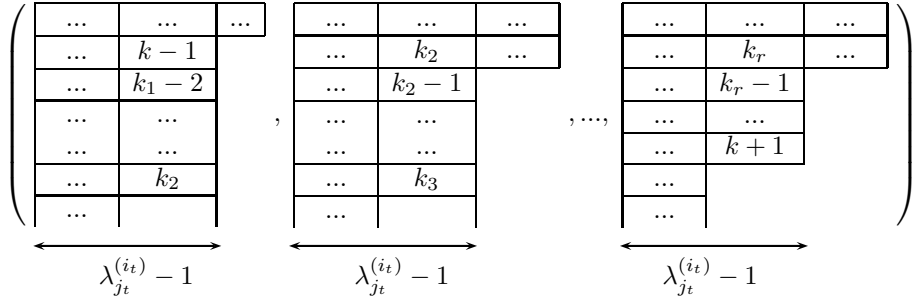
$$\left( \ \begin{array}{|c|c|c|}
                        \hline
                          \ \ ... \ \     & ... & ... \\  
                        \hline
                          \ \   ...\  \   &   k-1    \\ 
                      \cline{1-2}
                         \ \  ...\ \      & k_1-2 \\        
                         \cline{1-2}
              \ \   ...  \ \ & ... \\  
              \ \ ... \ \  & \ \ ... \ \   \\ 
                      \cline{1-2}
                         \ \  ...\ \       & k_2 \\ 
                        \cline{1-2}
               \ \   ...\ \      &  \\
               
                           \end{array}\ 
,\  \begin{array}{|c|c|c|}
                        \hline
                          \ \ ... \ \     & ... & ... \\  
                        \hline
                          \ \   ...\  \   &   k_2   &\ \  ... \ \ \\ 
                      \hline
                         \ \  ...\ \      & k_2-1 \\        
                         \cline{1-2}
              \ \   ...  \ \ & ... \\  
              \ \ ... \ \  & \ \ ... \ \   \\ 
                      \cline{1-2}
                         \ \  ...\ \       & k_3 \\ 
                        \cline{1-2}
               \ \   ...\ \      &  \\
               
                           \end{array}\ 
,...,\  \begin{array}{|c|c|c|}
                        \hline
                          \ \ ... \ \     & ... & ... \\  
                        \hline
                          \ \   ...\  \   &   k_r   &\ \  ... \ \ \\ 
                      \hline
                         \ \  ...\ \      & k_r-1 \\        
                         \cline{1-2}
              \ \   ...  \ \ & ... \\ 
              \cline{1-2}
              \ \ ... \ \  & k+1  \\ 
                      \cline{1-2}
                         \ \  ...\ \  \\ 
                        \cline{1-1}
               \ \   ...\ \        \\
               
                           \end{array}\ 
 \right)  $$
\begin{pspicture}(20,0.4)
\psline{<->}(0.4,0.4)(2.9,0.4)
\psline{<->}(3.8,0.4)(6.3,0.4)
\psline{<->}(8.2,0.4)(10.7,0.4)
\rput(1.9,0){$\lambda_{j_t}^{(i_t)}-1$}
\rput(5.3,0){$\lambda_{j_t}^{(i_t)}-1$}
\rput(9.4,0){$\lambda_{j_t}^{(i_t)}-1$}
\end{pspicture}\\
\\
Now, we can delete the remaining removable $k$-nodes of 
length $\lambda_{j_t}^{(i_t)}$ without violating the second condition.
Repeating the same argument for those $\xi_u$ with 
$\lambda_{j_u}^{(i_u)}>\lambda_{p_1}^{(l_1)}$, we conclude 
that the second condition holds.\\
Thus,   $\underline{\lambda}'$ is a FLOTW $d$-partition.
\end{proof}

Thanks to this lemma, we can now associate to each FLOTW 
$d$-partition  a residue sequence which have ``good'' properties
 according to  the $a$-function.

 Let $\underline{\lambda}$ be a  FLOTW $d$-partition associated 
to $\{e;v_0,...,v_{d-1}\}$.
 By using  Lemma  \ref{deuxiemelemme}, there exists a removable 
node $\xi_1$ with residue  $k$ on a part  $\lambda_{j_1}^{(i_1)}$ 
with maximal length, such that there doesn't exist a $k-1$-node 
on the border of a part with the same length as the length of the 
part of  $\xi_1$.

Let  $\gamma_1$, $\gamma_2$,..., $\gamma_r$ be the  $k-1$-nodes on 
the border of   $\underline{\lambda}$, on parts  
$\displaystyle{\lambda_{p_1}^{(l_1)}}$ ${\displaystyle{\geq
{\lambda_{p_2}^{(l_2)}}}}$ $\displaystyle{\geq{...}}$ 
$\displaystyle\geq{\lambda_{p_r}^{(l_r)}}$. Then,  Lemma 
\ref{premierlemme} implies that:
$$\lambda_{j_1}^{(i_1)}-j_1+m^{(i_1)}-1>\lambda^{(l_k)}_{p_k}-p_k+m^{(l_k)}
\qquad{k=1,...,r}.$$

As in  Lemma \ref{troisiemelemme}, let   $\xi_1$, $\xi_2$,..., $\xi_{u}$
 be the  removable $k$-nodes  on the border of   $\underline{\lambda}$
 on parts  $\displaystyle{\lambda_{j_1}^{(i_1)}\geq{\lambda_{j_2}^{(i_2)}
\geq{...}\geq{\lambda_{j_u}^{(i_u)}}}}$   such that:
$$\lambda_{j_u}^{(i_u)}>\lambda_{p_1}^{(l_1)}.$$
By  Lemma \ref{premierlemme}, we have: 
$$ \lambda_{j_t}^{(i_t)}-j_t+m^{(i_t)}-1>\lambda^{(l_k)}_{p_k}-p_k
+m^{(l_k)} \qquad{t=1,...,u}\qquad{k=1,...,r}.$$

We remove the nodes  $\xi_1$, $\xi_2$,..., $\xi_u$ from
 $\underline{\lambda}$. Let  $\underline{\lambda}'$ be the 
 resulting $d$-partition.   $\underline{\lambda}'$ is a FLOTW 
 $d$-partition by Lemma \ref{troisiemelemme}.

Now  the  $a-\sequence$  of residues of  $\underline{\lambda}$
 is defined recursively as follows:
\begin{df}\label{asequence}
   $a-\sequence (\underline{\lambda})=a-\sequence 
(\underline{\lambda}'),\underbrace{k,...,k}_{u}.$
\end{df}

Note that if we have 
 $a-\sequence (\underline{\lambda}')=a-\sequence 
(\underline{\lambda}''),\underbrace{k',...,k'}_{u'}$
for some $k'\in{\{0,1,...,e-1\}}$ and $u'\in{\mathbb{N}}$,
 then we have $k\neq{k'}$. Indeed, let $\xi_1 '$ be a removable
 $k'$-node on a part  ${\lambda'}_{a_1}^{(b_1)} $ with maximal 
length, such that there doesn't exist a $k'-1$-node on the border
 of a part of length  ${\lambda'}_{a_1}^{(b_1)}$. Assume that $k=k'$. 
We can't have  ${\lambda'}_{a_1}^{(b_1)}=\lambda_{j_1}^{(i_1)}$, 
otherwise $\xi_1'=\xi_t$ for some $t\in{\{1,2,...,u\}}$, a contradiction.
 Thus, we have  ${\lambda'}_{a_1}^{(b_1)}<\lambda_{j_1}^{(i_1)}$ 
and as  $\xi_1'\neq\xi_t$ for all $t\in{\{1,2,...,u\}}$, there exists
 a  $k'-1$-node on the border of a part of length  
${\lambda'}_{a_1}^{(b_1)}$ in $\underline{\lambda}'$, a contradiction
 again.
\begin{exa} 
We assume that $d=2$, $v_0=0$, $v_1=1$ and   $e=4$. The FLOTW 
 $2$-partitions are the  $2$-partitions $(\lambda^{(0)},\lambda^{(1)})$
 which satisfy:
\begin{itemize}
\item for all $i\in{\mathbb{N}_{>0}}$:
$$\lambda_i^{(0)}\geq{\lambda^{(1)}_{i+1}},$$
$$\lambda^{(1)}_i\geq{\lambda^{(0)}_{i+3}};$$
\item  for all  $k>0$, among the residues appearing at the right 
ends of the length $k$ rows of   $\underline\lambda$, at least one 
element of  $\{0,1,2,3\}$ does not occur.
\end{itemize}

We consider the  $2$-partition $\underline{\lambda}=(2.2,2.2.1)$ 
with the following diagram:
$$  \left( \ \begin{tabular}{|c|c|}
                         \hline
                           0  & 1  \\
                          \hline
                           3 & 0 \\
                          \hline
                             \end{tabular}\  ,\
                       \begin{tabular}{|c|c|}
                         \hline
                           1  & 2   \\
                          \hline
                          0   & 1  \\
                          \hline
                        3\\
                        \cline{1-1}
                             \end{tabular}\
                                                 \right)$$
\\
$\underline{\lambda}$ is a FLOTW $2$-partition.

We search the $a-\sequence$ of $\underline{\lambda}$: we have to find 
 $k\in{\{0,1,2,3\}}$, $u\in{\mathbb{N}_{>0}}$ and a $2$-partition 
$\underline{\lambda}'$ such that:
$$a-\sequence(\underline{\lambda})=a-\sequence(\underline{\lambda}'),
\underbrace{k,...,k}_{u}.$$
The parts with maximal length are the parts with length  $2$ and 
the residues of the removable nodes on these parts are $1$ and $0$.

For $k=1$, we have a part of length $2$  with a node on the border 
which have  $k-1=0\ (\textrm{mod}\ e)$ as  a residue.

So, we have to take  $k=0$  and we can remove this node as 
$3$-node(s) on the border have maximal length $1$.
There is no other removable  $0$-node in  $\underline{\lambda}$, so:
$$a-\sequence (\underline{\lambda})=a-\sequence(2.1,2.2.1),0.$$
We can verify that the  $2$-partition $(2.1,2.2.1)$ is a FLOTW  $2$-partition.

Now, the removable nodes on the part with maximal length have $1$ 
as a residue and there is no $0$-node, so:
$$a-\sequence(\underline{\lambda})=a-\sequence(1.1,2.1.1),1,1,0.$$
Repeating the same procedure, we obtain:
$$a-\sequence(\underline{\lambda})=1,0,0,3,3,2,1,1,0.$$
\end{exa}

 Let   $\underline{\lambda}$ be a  $d$-composition of rank  $n$ 
and let $k$ be a residue. Following the notations of paragraphe 
\ref{not}, let    $\underline{\mu}$ be a  $d$-composition of rank
  $n+1$ such that:

 $$\underline{\lambda}{\overset{k}{\underset{(j,p)}{\mapsto}}}
{\underline{\mu}}.$$
Then, we write:
 $$\underline{\lambda}{\overset{k\textrm{-opt}}
{\underset{(j,p)}{\mapsto}}}{\underline{\mu}},$$
if we have:
$$\lambda^{(p')}_{j'}-j'+m^{(p')}\leq{\lambda^{(p)}_{j}-j+m^{(p)}},$$
 for all $d$-composition $\underline{\mu'}$ of rank $n+1$ such that 
 $\underline{\lambda}{\overset{k}{\underset{(j',p')}{\mapsto}}}
{\underline{\mu}'}$. 
\begin{re}\label{partition}
Let $k\in{\{0,1,...,e-1\}}$. Assume that $\underline{\lambda}$
 is a $d$-partition and that  $\underline{\mu}$ and  $\underline{\mu}'$
 are $d$-compositions such that:
 $$\underline{\lambda}{\overset{k\textrm{-opt}}{\underset{(j,p)}
{\mapsto}}}{\underline{\mu}}\qquad{\textrm{and}}\qquad\underline{\lambda}
{\overset{k\textrm{-opt}}{\underset{(j',p')}{\mapsto}}}{\underline{\mu}'}$$
This implies that there exists $t\in{\mathbb{Z}}$ such that:
$$\lambda^{(p')}_{j'}-j'+v_{p'}={\lambda^{(p)}_{j}-j+v_{p}+te}
\equiv k\ (\textrm{mod}\ e).$$
Moreover, we have:
$$\lambda^{(p')}_{j'}-j'+m^{(p')}={\lambda^{(p)}_{j}-j+m^{(p)}}.$$
We obtain:
$$te=(p'-p)\frac{e}{d}.$$
Hence, we have $p=p'$ and since  $\underline{\lambda}$ is a $d$-partition, 
we have $(j,p)=(j',p')$ and ${\underline{\mu}'}={\underline{\mu}}$.

\end{re}

\begin{prp}\label{agr} 
Let  $\underline{\lambda}$ be a FLOTW   $d$-partition of rank $n$ and let:
 $$ a-\sequence (\underline{\lambda})=  s_1,s_2,...,s_n.$$
  Then, there exists a sequence of FLOTW $d$-partitions
 $\underline{\lambda}[0]=\underline{\emptyset}$,
 $\underline{\lambda}[1]$,...,  $\underline{\lambda}[n]=\underline{\lambda}$
 such that for all  $l\in{\{0,...,n-1\}}$, we have:
 $$\underline{\lambda}[l]{\overset{s_l\textrm{-opt}}
{\underset{(j_l,p_l)}{\mapsto}}}{\underline{\lambda}[l+1]}.$$
 
\end{prp}
\begin{proof} 
Let $\underline{\lambda}$ be a FLOTW $d$-partition. Assume that:
  $$\displaystyle a-\sequence (\underline{\lambda})=a-\sequence 
(\underline{\lambda}')\underbrace{i_s,...,i_s}_{a_s},$$
where $\underline{\lambda}'$ is the FLOTW $d$-partition as in 
Definition \ref{asequence}. Define $\underline\lambda [n-a_s] := 
\underline{\lambda}'$. Then, there exist $d$-compositions   
$\underline\lambda [i]$, $i=n-a_s+1,...,n$ such that
$$\underline{\lambda}[n-a_s]{\overset{i_s\textrm{-opt}}{\mapsto}}
\underline{\lambda}[n-a_s+1]{\overset{i_s\textrm{-opt}}{\mapsto}} 
...{\overset{i_s\textrm{-opt}}{\mapsto}}\underline{\lambda}[n-1]
{\overset{i_s\textrm{-opt}}{\mapsto}}{\underline{\lambda}[n]}.$$
(we omit the indices of the nodes here to simplify the notations.)

We claim that $\underline{\lambda}[n]=\underline{\lambda}$. Indeed, 
by the discussion above Definition \ref{asequence}, the nodes of 
$\underline{\lambda}[n]/\underline{\lambda}[n-a_s]$ are precisely 
the nodes of  $\underline{\lambda}/\underline{\lambda}'$. Note also 
that in the above graph, the nodes are succesively added on the 
greatest part where we have an addable $i_s$-node. Moreover, as 
$\underline{\lambda}$ is a $d$-partition, it follows that for all 
$i\in{\{n-a_s,n-a_s+1,...,n\}}$, $\underline{\lambda}[i]$ is a $d$-partition. 

Continuing in this way,  there exist  $d$-partitions  
$\underline{\lambda}[l]$ ($l=1,...,n$) such that:
$$\underline{\emptyset}{\overset{i_1\textrm{-opt}}{\mapsto}}
\underline{\lambda}[1]{\overset{i_1\textrm{-opt}}{\mapsto}} 
...{\overset{i_1\textrm{-opt}}{\mapsto}}\underline{\lambda}[a_1]
{\overset{i_2\textrm{-opt}}{\mapsto}} ...{\overset{i_2\textrm{-opt}}
{\mapsto}\underline{\lambda}[a_1+a_2] {\overset{i_3-\textrm{opt}}
{\mapsto}}...{\overset{i_s\textrm{-opt}}{\mapsto}}\underline{\lambda}[n]}
\ \ \ (1).$$

 By Lemma \ref{troisiemelemme} and by definition of the  
$a-\sequence$, $\underline{\lambda}[a_1]$,  $\underline{\lambda}
[a_1+a_2]$,..., $\underline{\lambda}[n]$ are FLOTW  $d$-partitions. 
Now, assume that there exists a  $d$-partition  $\underline{\mu}$ 
which appears in  $(1)$ and which is not a FLOTW $d$-partition. 
Then, there exists $r\in{\{1,...,s\}}$ such that:
$$\underline{\lambda}[a_1+...+a_{r-1}]{\overset{i_r\textrm{-opt}}
{\mapsto}}...{\overset{i_r\textrm{-opt}}{\mapsto}}\ \underline{\mu}
\ {\overset{i_r\textrm{-opt}}{\mapsto}}...    
 {\overset{i_r\textrm{-opt}}{\mapsto}}\underline{\lambda}[a_1+...+a_{r}]
\ \ \ (2)$$
We denote  $\underline{\nu}:=   \underline{\lambda}[a_1+...+a_{r-1}]$. 
Then, we have two  cases to consider.
\begin{itemize}
\item If $\underline{\mu}$ violates the first condition to be a 
 FLOTW $d$-partition, then  there exist  $j\in{\mathbb{N}_{>0}}$
 and  $i\in{\{0,1,...,d-2\}}$ such that   $\mu^{(i)}_{j}<
\mu^{(i+1)}_{j+v_{i+1}-v_i}$ (the case $\mu^{(d-1)}_{j}<
\mu^{(0)}_{j+e+v_{0}-v_{d-1}}$ is similar to this one). It 
implies that we have added an  $i_r$-node on  $\nu^{(i+1)}_{j+v_{i+1}-v_i}$
 and that $\nu^{(i)}_{j}=\nu^{(i+1)}_{j+v_{i+1}-v_i}$. We obtain:
 $$\nu^{(i)}_{j}-j+m^{(i)}>\nu^{(i+1)}_{j+v_{i+1}-v_i}-(j+v_{i+1}-v_i)
+m^{(i+1)}.$$
 It implies that we can't add any  $i_r$-node on   $\mu^{(i)}_{j}$
 in  $(2)$ because we add $k$-nodes by decreasing order with respect
 to $\nu^{(i')}_{j'}-j'+m^{(i')}$. Thus, $\underline{\lambda}[a_1+...+a_r]$
 is not a FLOTW  $d$-partition. This is a contradiction.
\item Assume that  $\underline{\mu}$ is the first place where the 
 second condition  to be a  FLOTW $d$-partition is violated. Thus,
 there exist $a_1+...+a_{r-1}<j<a_1+...+a_r$ such that
 $\underline{\mu}=\underline{\mu}[j]$  and $\underline{\mu}[j-1]$ 
is a FLOTW $d$-partition. In $(2)$, we have added an  $i_r$-node
 on a part   $\nu^{(i)}_j$  and   the set of residues appearing 
at the right ends of parts of length $\nu^{(i)}_j+1$  of  
 $\underline{\mu}$ is equal to $\{0,...,e-1\}$. If there is a
 $i_r-1$-node on the border of a part of length greater than 
$\nu_j^{(i)}$ which will be occupied before reaching 
$\underline{\lambda}[a_1+...+a_{r}]$, Lemma \ref{premierlemme}
 implies that this node had to be occupied before the $i_r$-node
 on the part   $\nu^{(i)}_j$ is added. Thus,  the set of residues 
appearing at the right ends of parts of length $\nu^{(i)}_j+1$  
in  $\underline{\lambda}[a_1+...+a_r]$ is equal to $\{0,...,e-1\}$. 
This is a contradiction as  $\underline{\lambda}[a_1+...+a_r]$ is 
a FLOTW $d$-partition.
\end{itemize}
\end{proof}

Now, we  associate a graph with a  FLOTW $d$-partition as follows:

\begin{df}\label{agraph} Let  $\underline{\lambda}$ be a 
FLOTW $d$-partition of rank $n$, let  $s_1,s_2,...,s_n$ be 
its $a-\sequence$  of residues. Then by Remark \ref{partition} 
and Proposition \ref{agr}, there exist unique $d$-partitions 
$\underline{\lambda}[l]$ ($l=1,...,n$) such that: 
$$\underline{\emptyset}{\overset{s_1\textrm{-opt}}
{\underset{(j_1,p_1)}{\mapsto}}\underline{\lambda}[1]} 
{\overset{s_2\textrm{-opt}}{\underset{(j_2,p_2)}{\mapsto}}
\underline{\lambda}[2]}...{\overset{s_n\textrm{-opt}}{\underset{(j_n,p_n)}
{\mapsto}}\underline{\lambda}{[n]}=\underline{\lambda}},$$
where $\underline{\emptyset}$ is the empty $d$-partition. We call 
this the $a$-graph of $\underline{\lambda}$.
\end{df}
\begin{exa} 
Assume that $d=2$, $v_0=0$, $v_1=1$  and   $e=4$. Let 
$\underline{\lambda}=(2.2,2.2.1)$.
$$a-\sequence (\underline{\lambda})=1,0,0,3,3,2,1,1,0.$$
Then, the  $a$-graph associated to   $\underline{\lambda}$ 
is the following one:

$$(\emptyset,\emptyset){\overset{1\textrm{-opt}}{\underset{(1,1)}{\mapsto}}}
(\emptyset,1){\overset{0\textrm{-opt}}{\underset{(1,0)}
{\mapsto}}}(1,1){\overset{0\textrm{-opt}}{\underset{(2,1)}
{\mapsto}}}(1,1.1){\overset{3\textrm{-opt}}{\underset{(2,0)}
{\mapsto}}}(1.1,1.1){\overset{3\textrm{-opt}}{\underset{(3,1)}
{\mapsto}}}(1.1,1.1.1)\overset{2\textrm{-opt}}{\underset{(1,1)}
{\mapsto}}$$
$$\overset{2\textrm{-opt}}{\underset{(1,1)}
{\mapsto}}(1.1,2.1.1)\overset{1\textrm{-opt}}{\underset{(1,0)}
{\mapsto}}(2.1,2.1.1)\overset{1\textrm{-opt}}{\underset{(2,1)}
{\mapsto}}(2.1,2.2.1)\overset{0\textrm{-opt}}{\underset{(2,0)}
{\mapsto}}(2.2,2.2.1)$$
\end{exa}
Now, we deal with  a fundamental property concerning the $a$-graph o
f a FLOTW $d$-partition. First, we have to show some properties 
about $a$-functions.

\begin{df} 
Let  $\underline{\mu}$ and  $\underline{\nu}$ be  $d$-compositions 
of rank  $n$. Let   $\textbf{B}_{\underline{\mu}}$ and 
 $\textbf{B}_{\underline{\nu}}$  be two ordinary symbols of 
 $\underline{\mu}$ et  $\underline{\nu}$ with the same height. 
Let  $\textbf{B}_{ \underline{\mu}}[m]'$,
  $\textbf{B}_{ \underline{\nu}}[m]'$ be as in paragraph 
\ref{symbole}. Then we write:
$$\underline{\mu} \prec{\underline{\nu}},$$
if:
$$\sum_{{0\leq{i}\leq{j}<d}\atop{{(a,b)\in{B_{\underline{\mu}}'^{(i)}
\times{B_{\underline{\mu}}'^{(j)}}}}\atop{a>b\ \textrm{if}\ i=j}}}
{\min{\{a,b\}}}   -\sum_{{0\leq{i,j}<d}
\atop{{a\in{B_{\underline{\mu}}'^{(i)}}}
\atop{1\leq{k}\leq{a}}}}{\min{\{k,m^{(j)}\}}}<$$
$$\sum_{{0\leq{i}\leq{j}<d}\atop{{(a,b)\in{B_{\underline{\nu}}'^{(i)}
\times{B_{\underline{\nu}}'^{(j)}}}}\atop{a>b\ \textrm{if}\ i=j}}}
{\min{\{a,b\}}}   -\sum_{{0\leq{i,j}<d}\atop{{a\in{B_{\underline{\nu}}
'^{(i)}}}\atop{1\leq{k}\leq{a}}}}{\min{\{k,m^{(j)}\}}}.$$
In particular, if  $\underline{\mu}$ and   $\underline{\nu}$ are some 
 $d$-partitions, we have:
$$\underline{\mu} \prec{\underline{\nu}}\iff{a(\underline{\mu})
<a(\underline{\nu})}.$$
\end{df}

The following results give  the consequences of the properties 
showed in this section in terms of $a$-function:

\begin{lmm}\label{propaf} Let  $\underline{\lambda}$ be a 
 $d$-composition of rank $n$, let  $\textbf{B}:=(B^{(1)},...,B^{(l-1)})$ 
be an ordinary symbol of $\underline{\lambda}$,  let  $\beta_1$ 
and $\beta_2$ be two elements of   $\textbf{B}[m]'$, we assume
 that:
$$\beta_1 < \beta_2.$$
 Let  $l\in{\mathbb{N}_{>0}}$. We add  $l$ nodes  to $\underline{\lambda}$
 on the part associated to    $\beta_1$. Let  $\underline{\mu}$
 be the resulting $d$-composition of rank $n+l$.  We add  $l$ nodes
 to $\underline{\lambda}$  on the part associated to    $\beta_2$.
 Let  $\underline{\nu}$ be the resulting $d$-composition of rank 
$n+l$.  Let  $\textbf{B}_{\underline{\mu}}$ and 
  $\textbf{B}_{\underline{\nu}}$ be the two ordinary symbols
 of  $\underline{\mu}$ and   $\underline{\nu}$ with the same 
height as   $\textbf{B}$ and let  $\textbf{B}_{ \underline{\mu}}[m]'$, 
 $\textbf{B}_{ \underline{\nu}}[m]'$ be as in the  paragraph 
\ref{symbole}.   Then, we have:
$$\underline{\nu} \prec{\underline{\mu}}.$$
In particular, if   $\underline{\mu}$ and   $\underline{\nu}$ 
are some  $d$-partitions, then:
$$a(\underline{\mu})>a(\underline{\nu}).$$
\end{lmm}
\begin{proof} 
 We assume that $\beta_1$ is on   $B^{(i_1)}$ and  that $\beta_2$ 
is on   $B^{(i_2)}$. Then, we have:

$$\sum_{{0\leq{i}\leq{j}<d}\atop{{(a,b)\in{B_{\underline{\mu}}'^{(i)}
\times{B_{\underline{\mu}}'^{(j)}}}}\atop{a>b\ \textrm{if}\ i=j}}}
{\min{\{a,b\}}}   -\sum_{{0\leq{i,j}<d}
\atop{{a\in{B_{\underline{\mu}}'^{(i)}}}
\atop{1\leq{t}\leq{a}}}}{\min{\{t,m^{(j)}\}}}$$
$$-(\sum_{{0\leq{i}\leq{j}<d}\atop{{(a,b)\in{B'^{(i)}
\times{B}'^{(j)}}}\atop{a>b\ \textrm{if}\ i=j}}}{\min{\{a,b\}}}
  -{}\sum_{{0\leq{i,j}<d}\atop{{a\in{B'^{(i)}}}\atop{1\leq{t}\leq{a}}}}
{\min{\{t,m^{(j)}\}}})=$$

$$=\sum_{{a\in{B^{(i)}}\atop{0\leq{i}\leq{d-1}}}\atop{a\neq{\beta_1}
\ \textrm{if}\ i=i_1}}{(\min{\{a,\beta_1+l\}}-\min{\{a,\beta_1\}})}
{}-\sum_{{0\leq{j}\leq{d-1}}\atop{\beta_1<t\leq{\beta_1+l}}}
{\min{\{t,m^{(j)}\}}}.$$
and we have an analogous formula for  $\underline{\nu}$.\\
Now, for  $a\in{B^{(j)}}$,  $j=0,...,d-1$, we have:

$$\min{\{a,\beta_1+l\}}-\min{\{a,\beta_1\}}=\left\lbrace
                                                        \begin{array}{ll}
        0 & \textrm{if }a<\beta_1,\\
     a-\beta_1 & \textrm{if }\beta_1\leq{a}<\beta_1+l,\\
   l   &  \textrm{if }{a}\geq{\beta_1+l}.\\
\end{array}\right.$$
As $\beta_1<\beta_2$, we have:
$$\sum_{{a\in{B^{(i)}}\atop{0\leq{i}\leq{d-1}}}\atop{a\neq{\beta_1}
\ \textrm{if}\ i=i_1}}{(\min{\{a,\beta_1+l\}}-\min{\{a,\beta_1\}})}>$$
$$\sum_{{a\in{B^{(i)}}\atop{0\leq{i}\leq{d-1}}}\atop{a\neq{\beta_2}
\ \textrm{if}\ i=i_2}}{(\min{\{a,\beta_2+l\}}-\min{\{a,\beta_2\}})}.$$
Moreover, $\min{\{t,m^{(j)}\}}\leq \min{\{t+\beta_2-\beta_1,m^{(j)} \}}$
 implies that:
$$\sum_{{0\leq{j}\leq{d-1}}\atop{\beta_1<t\leq{\beta_1+l}}}
{\min{\{t,m^{(j)}\}}}\leq{\sum_{{0\leq{j}\leq{d-1}}
\atop{\beta_2<t\leq{\beta_2+l}}}{\min{\{t,m^{(j)}\}}}}.$$
\end{proof}

\begin{prp}\label{propaf2}
Let  $\underline{\lambda}[n]$ be a FLOTW  $d$-partition  of rank 
 $n$ and  $s_1,s_2,...,s_n$ be its  $a-\sequence$  of residues. 
We consider the $a$-graph of  $\underline{\lambda}[n]$: 
$$(\emptyset,\emptyset){\overset{s_1\textrm{-opt}}
{\underset{(j_1,p_1)}{\mapsto}}\underline{\lambda}[1]} 
{\overset{s_2\textrm{-opt}}{\underset{(j_2,p_2)}{\mapsto}}
\underline{\lambda}[2]}...{\overset{s_n\textrm{-opt}}
{\underset{(j_n,p_n)}{\mapsto}}\underline{\lambda}[n]},$$
where all the $d$-partitions appearing in this graph are FLOTW 
$d$-partitions. Then, if we have another graph of  $d$-compositions:
$$(\emptyset,\emptyset){\overset{s_1}{\underset{(j'_1,p'_1)}{\mapsto}}
\underline{\mu}[1]} {\overset{s_2}{\underset{(j'_2,p'_2)}{\mapsto}}
\underline{\mu}[2]}...{\overset{s_n}{\underset{(j'_n,p'_n)}{\mapsto}}
{\underline{\mu}[n]}},$$
we have  $\underline{\lambda}[n]\prec{\underline{\mu}[n]}$ if 
$\underline{\lambda}[n]\neq{\underline{\mu}[n]}$.
\end{prp}
\begin{proof} 
We will argue by   induction on $n-r\in{\mathbb{N}}$. Assume that
 $\underline{\lambda}[r-1]=\underline{\mu}[r-1]$. We want to show 
 that if we have the  following graphs:
$$\underline{\lambda}[r-1]{\overset{s_r\textrm{-opt}}
{\underset{(j_r,p_r)}{\mapsto}}\underline{\lambda}[r]}
{\overset{s_{r+1}\textrm{-opt}}{\underset{(j_{r+1},p_{r+1})}
{\mapsto}}}   ...{\overset{s_n\textrm{-opt}}{\underset{(j_n,p_n)}
{\mapsto}}\underline{\lambda}[n]}\ \ \ \ (1)$$
$$\underline{\mu}[r-1]{\overset{s_r}{\underset{(j'_r,p'_r)}{\mapsto}}
\underline{\mu}[r]} {\overset{s_{r+1}}{\underset{(j'_{r+1},p'_{r+1})}
{\mapsto}}}...{\overset{s_n}{\underset{(j'_n,p'_n)}{\mapsto}}
\underline{\mu}[n]}\ \ \ \ (2)$$
then $\underline{\lambda}[n]\prec{\underline{\mu}[n]}$ if 
$\underline{\lambda}[n]\neq{\underline{\mu}[n]}$.
\\
\underline{Assume that $n-r=0$}:\\
Assume that $\underline{\lambda}[n]\neq{\underline{\mu}[n]}$.
 To simplify the notations, we write $\underline{\lambda}:=
\underline{\lambda}[n-1]=\underline{\mu}[n-1]$. We have:
$$   \underline{\lambda}{\overset{s_n\textrm{-opt}}
{\underset{(j_n,p_n)}{\mapsto}}\underline{\lambda}[n]}
\qquad{\textrm{and}}\qquad   \underline{\lambda}{\overset{s_n}
{\underset{(j'_n,p'_n)}{\mapsto}}\underline{\mu}[n]}$$
As $\underline{{\lambda}}$ is a $d$-partition, by Remark 
\ref{partition}, we have
$$\lambda^{(p)}_{j}-j+m^{(p)}>{\lambda^{(p')}_{j'}-j'+m^{(p')}}$$
Let  $\textbf{B}_{\underline{\lambda}}$ be  an ordinary symbol 
associated to $\underline{\lambda}$. Then, $\underline{\lambda}[n]$ 
and  ${\underline{\mu}[n]}$ are obtained from  $\underline{\lambda}$ 
by adding a node on parts  associated to two elements $\beta_1$ and 
 $\beta_2$  of $\textbf{B}_{\underline{\lambda}}[m]'$. We have:
$$\beta_1>\beta_2.$$  
 We are in the setting of Lemma \ref{propaf}, hence:
 $$\underline{\lambda}[n]\prec{\underline{\mu}[n]}.$$
\\
\underline{Assume that $n-r>0$}:\\
 Let $t$ be such that $r-1\leq t<n$. If the residues of the
 right ends of the parts $(j_r,p_r)$ and  $(j'_r,p'_r )$ of
 $\underline{\mu}[t]$ are the same, we say that $t$ is admissible.
 Let $t<t'$ be the two first consecutive admissible indices. 
Suppose that the lengths of the part $(j_r,p_r)$ (resp. $(j'_r,p'_r)$) 
increases by $N$ (resp. $N'$) between $\underline{\mu}[t]$ 
and  $\underline{\mu}[t']$. If the first node added to 
 $\underline{\mu}[t]$ is on the $(j_r,p_r)$-part, we do nothing.
 Otherwise, we add $N$ nodes to the  $(j'_r,p'_r)$-part and
 $N'$ nodes to the  $(j_r,p_r)$-part instead of adding  $N$ nodes
 to the  $(j_r,p_r)$-part and $N'$ nodes to the  $(j'_r,p'_r)$-part.

Next we consider the consecutive admissible indices  $t'<t''$ and 
repeat the same procedure until we reach the final consecutive admissible indices. Then we get a new graph as follows:
$$\underline{\lambda}[r-1]{\overset{s_r\textrm{-opt}}
{\underset{(j_r,p_r)}{\mapsto}}\underline{\lambda}[r]} 
{\overset{s_{r+1}}{\underset{(j''_{r+1},p''_{r+1})}
{\mapsto}}}...{\overset{s_n}{\underset{(j''_n,p''_n)}
{\mapsto}}\underline{\nu}[n]},\ \ \ \ (3)$$
In  $\underline{\nu}[n]$ and  $\underline{\mu}[n]$,  only 
two parts ($(j_r,p_r)${}-part  and   $(j'_r,p'_r)$-part) may differ.
\begin{itemize}
\item If $(j_r,p_r)=(j'_r,p'_r)$, we conclude by induction.
\item If $(j_r,p_r)\neq{(j'_r,p'_r)}$,  to simplify, we write 
 $\underline{\lambda}:= \underline{\lambda}[n]$,  $\underline{\mu}
:= \underline{\mu}[n]$ and $\underline{\nu}:= \underline{\nu}[n]$.
Assume that   $\underline{\nu}\neq{\underline{\mu}}$, then we are
 in the setting of Lemma  \ref{propaf}: let  $\underline{\gamma}=
(\gamma^{(0)},...,\gamma^{(l-1)})$ be the  $d$-composition defined by:
\begin{align*}
&\gamma_j^{(p)}:=\mu^{(p)}_j & \textrm{if }(j,p)\neq (j'_r,p'_r),\\
&\gamma_{j'_r}^{(p'_r)}:=\nu^{(p'_r)}_{j'_r}.&
\end{align*}
Then, $\underline{\mu}$ is obtained from $\underline{\gamma}$ 
by adding nodes on   $\gamma_{j'_r}^{(p'_r)}$ where as $\underline{\nu}$ 
is obtained from $\underline{\gamma}$ by adding the same number of 
nodes on   $\gamma_{j_r}^{(p_r)}$. Now, we know that:
$$\lambda_{j_r}^{(p_r)}-j_r+m^{(p_r)}>\lambda_{j'_r}^{(p'_r)}-
j'_r+m^{(p'_r)}.$$ This implies: 
$$\gamma_{j_r}^{(p_r)}-j_r+m^{(p_r)}>\gamma_{j'_r}^{(p'_r)}-
j'_r+m^{(p'_r)}.$$
This follows from the fact that  the residue associated to $(j_r,p_r)$ 
and $(j'_r,p'_r)$ in $\underline{\lambda}$ are the same.

Thus,    $\underline{\nu}$ and  ${\underline{\mu}}$ are obtained from 
 $\underline{\gamma}$ by adding the same number of nodes on parts  
associated to two elements $\beta_1$ and  $\beta_2$  of 
$\textbf{B}_{\gamma}[m]'$ where $\textbf{B}_{\gamma}$ is an
 ordinary symbol associated to $\gamma$. The above discussion 
shows that:
 $$\beta_1>\beta_2.$$
So, by Lemma \ref{propaf}, we have:
$$\underline{\nu}\prec{\underline{\mu}}.$$
By induction hypothesis, we have :
$$\underline{\lambda}\prec{\underline{\nu}}\qquad{\textrm{or}}
\qquad \underline{\lambda}={\underline{\nu}}.$$
We conclude that:
$$\underline{\lambda}\prec{\underline{\mu}}.$$
\end{itemize}
\end{proof}
The following proposition is now clear.
\begin{prp}\label{resu}
Let $\underline{\lambda}[n]$ be a  FLOTW $d$-partition of rank 
 $n$ and  let  $s_1,s_2,...,s_n$ be its  $a-\sequence$ of residues. 
We consider the $a$-graph of  $\underline{\lambda}[n]$: 
$$(\emptyset,\emptyset){\overset{s_1\textrm{-opt}}
{\underset{(j_1,p_1)}{\mapsto}}\underline{\lambda}[1]} 
{\overset{s_2\textrm{-opt}}{\underset{(j_2,p_2)}{\mapsto}}
\underline{\lambda}[2]}...{\overset{s_n\textrm{-opt}}{\underset{(j_n,p_n)}
{\mapsto}}\underline{\lambda}[n]}$$
then, if we have another graph with  $d$-partitions:
$$(\emptyset,\emptyset){\overset{s_1}{\underset{(j'_1,p'_1)}{\mapsto}}
\underline{\mu}[1]} {\overset{s_2}{\underset{(j'_2,p'_2)}{\mapsto}}
\underline{\mu}[2]}...{\overset{s_n}{\underset{(j'_n,p'_n)}{\mapsto}}
{\underline{\mu}[n]}}$$
we have  $a(\underline{\lambda}[n])<a({\underline{\mu}[n]})$ if 
 $\underline{\lambda}[n]\neq{\underline{\mu}[n]}$.
\end{prp}

In the following paragraph, we see  consequences of this result: 
we give an interpretation of  the parametrization of the simple 
modules by Foda et al. for Ariki-Koike algebras in terms of 
$a$-functions.

\subsection{Principal result}
First, recall the   FLOTW order given in the paragraph \ref{FLOTW}, 
we have the following proposition.
\begin{prp}\label{orderFLOTW}   Let  $\underline{\lambda}$ be a 
FLOTW $d$-partition and let $\gamma=(a,b,c)$, $\gamma'=(a',b',c')$ 
be two nodes of  $\underline{\lambda}$ with the same residue. Then:
$$b-a+m^{(c)}>b'-a'+m^{(c')}\iff{\gamma \textrm{ is below } 
\gamma'\textrm{ with respect to the FLOTW order.}}$$
\end{prp}
\begin{proof} 
As $\gamma$ and $\gamma'$ have the same residue, there exists
 $t\in{\mathbb{Z}}$ such that:
$$b-a+v_c=b'-a'+v_{c'}+te.$$
Then, we have:
$$b-a+m^{(c)}>b'-a'+m^{(c')}\iff{(c-c')\frac{e}{d}<te}.$$
Then $t\geq{0}$ and if $t=0$, we have $c'>c$. So, $\gamma$ 
is below $\gamma'$.
\end{proof}

Now, we consider  the quantum group $\mathcal{U}_q$  of 
type $A^{(1)}_{e-1}$ as in the paragraph \ref{combinatoire}. 
Let $\overline{\mathcal{M}}$ be the   $\mathcal{U}_q$-module
 generated by the empty $d$-partition with respect to the JMMO
 action (see paragraph \ref{FLOTW}). 

First, we have the following result which gives the action of 
the divided powers on the multipartitions with respect to the 
JMMO  action. The demonstration is  analogous to the case of 
the AM action (see for example \cite[Lemma 6.15]{mathas} for $d=1$),
 we give the proof for the convenience of the reader.

\begin{prp}\label{mathas}
Let $\underline{\lambda}$ be a $d$-partition, $i\in{\{0,...,e-1\}}$ 
and $j$ a positive integer. Then:
$$f^{(j)}_i\underline{\lambda}=\sum{q^{\overline{N}_i^{b}
(\underline{\lambda},\underline{\mu})}\underline{\mu}},$$
where the sum is taken over all the $\underline{\mu}$ which
 satisfy:
$$\underline{\lambda}\underbrace{\overset{i}{\mapsto}...\overset{i}
{\mapsto}}_j\underline{\mu},$$
and where:
$$\overline{N}_i^{b}(\underline{\lambda},\underline{\mu}):=
\sum_{\gamma\in{[\underline{\mu}]/[\underline{\lambda}]}}\left(
 \sharp \left\{\begin{array}{c} \textrm{addable}\ i-\textrm{nodes}
\\\textrm{of}\ \underline{\mu}\ \textrm{below}\ \gamma\end{array}
 \right\}- \sharp \left\{\begin{array}{c}\textrm{removable}\ 
i-\textrm{nodes}\\\textrm{of}\  \underline{\lambda}\ \textrm{below}
\ \gamma\end{array}\right\}\right)$$
\end{prp}
\begin{proof}
The proof is by induction on $j\in{\mathbb{N}_{>0}}$:
\begin{itemize}
\item If $j=1$, this is the definition of the JMMO action 
(see paragraph \ref{FLOTW}).
\item If $j>1$, let $\underline{\mu}$ be a $d$-partition such that:
$$\underline{\lambda}\underbrace{{\overset{i}{\underset{(r_1,p_1)}
{\mapsto}}...{\overset{i}{\underset{(r_j,p_j)}{\mapsto}}}}}_j
\underline{\mu}.$$
We assume that the nodes $\gamma_l$ associated to $(r_l,p_l)$, 
$l=1,...,j$, are such that for all  $l=1,...,j-1$,  $\gamma_l$ 
is below $\gamma_{l+1}$. Let $\underline{\mu}_s$, $s=1,...,j$, 
be the $d$-partitions such that $[\underline{\mu}_s] / 
[\underline{\mu}_{s-1}]=\gamma_s$.

By induction, the coefficient of  $\underline{\mu}_s$ in 
$f^{(j-1)}_i\underline{\lambda}$ is $q^{\overline{N}_i^{b}
(\underline{\lambda},\underline{\mu}_s)}$. Hence, the 
coefficient of $\underline{\mu}$ in $f_if^{(j-1)}_i\underline{\lambda}$
 is $\displaystyle{\sum_{s=1}^jq^{\overline{N}_i^{b}
(\underline{\lambda},\underline{\mu}_s)+\overline{N}_i^{b}
(\underline{\mu}_s,\underline{\mu})}}$.

Now, we define the following numbers:
\end{itemize}
$$\overline{N}_i^{b}(\underline{\lambda},\underline{\mu},\gamma_s)
:=\left( \sharp \left\{\begin{array}{c} \textrm{addable}\ 
i-\textrm{nodes}\\\textrm{of}\ \underline{\mu}\ \textrm{below}\
 \gamma_s\end{array} \right\}- \sharp \left\{\begin{array}{c}
\textrm{removable}\ i-\textrm{nodes}\\\textrm{of}\ 
 \underline{\lambda}\ \textrm{below}\ \gamma_s\end{array}
\right\}\right)$$
\begin{itemize}
\item[]
where $s=1,...,j$.\\
We have:

\begin{eqnarray*} \overline{N}_i^{b}(\underline{\lambda}
,\underline{\mu}_s)& = &\sum_{t=1}^{s-1}\overline{N}_i^{b}
(\underline{\lambda},\underline{\mu}_s,\gamma_t)+\sum_{t=s+1}^{j}
 \overline{N}_i^{b}(\underline{\lambda},\underline{\mu}_s,\gamma_t) \\
 &=&  \sum_{t=1}^{s-1}\overline{N}_i^{b}(\underline{\lambda},
\underline{\mu},\gamma_t)+\sum_{t=s+1}^{j}( \overline{N}_i^{b}
(\underline{\lambda},\underline{\mu},\gamma_t)+1) \\
 &=&  \overline{N}_i^{b}(\underline{\lambda},\underline{\mu})-
\overline{N}_i^{b}(\underline{\lambda},\underline{\mu}_s,\gamma_s)+j-s.
\end{eqnarray*}
And:
$$\overline{N}_i^{b}(\underline{\mu}_s,\underline{\mu})=
\overline{N}_i^{b}(\underline{\lambda},\underline{\mu}_s,\gamma_s)
{}-s+1.$$

Hence, the coefficient of $\underline{\mu}$ in $f_if^{(j-1)}_i
\underline{\lambda}$ is $\displaystyle{q^{ \overline{N}_i^{b}
(\underline{\lambda},\underline{\mu})} \sum_{s=1}^jq^{j+1-2s}}$.

Now, we have $f_if^{(j-1)}_i\underline{\lambda}=[j]_qf^{(j)}_i$. 
Thus, the coefficient of   $\underline{\mu}$ in 
$f^{(j)}_i\underline{\lambda}$ is $q^{ \overline{N}_i^{b}
(\underline{\lambda},\underline{\mu})}$.
\end{itemize}
\end{proof}
Now, we  have the following result:
\begin{prp}\label{asuite} Let $\underline{\lambda}$ be a FLOTW 
$d$-partition and let $a-\sequence (\underline{\lambda})=
\underbrace{i_1,...,i_1}_{a_1},\underbrace{i_{2},...,i_{2}}_{a_{2}},
...,\underbrace{i_s,...,i_s}_{a_s}$ be its $a-\sequence$ of residues
 where for all $j=1,...,s-1$, we have $i_{j}\neq{i_{j+1}}$. Then, 
we have:
$$f^{(a_s)}_{i_s}f^{(a_{s-1})}_{i_{s-1}} ...f^{(a_1)}_{i_1}
\underline{\emptyset}=\underline{\lambda}+ \sum_{a(\underline{\mu})>
a(\underline{\lambda})}{c_{\underline{\lambda},\underline{\mu}}(q)
\underline{\mu}},$$
 where $c_{\underline{\lambda},\underline{\mu}}(q)\in{
\mathbb{Z}[q,q^{-1}]}$.
\end{prp}
\begin{proof}
By  Proposition \ref{resu}, we have:
$$f^{(a_s)}_{i_s}f^{(a_{s-1})}_{i_{s-1}} ...f^{(a_1)}_{i_1}
\underline{\emptyset}=c_{\underline{\lambda},\underline{\lambda}}(q) 
 \underline{\lambda}+ \sum_{a(\underline{\mu})>
a(\underline{\lambda})}{c_{\underline{\lambda},\underline{\mu}}
(q)\underline{\mu}},$$
where  $c_{\underline{\lambda},\underline{\mu}}
(q)\in{\mathbb{Z}[q,q^{-1}]}$. So, we have to show that:
 $$c_{\underline{\lambda},\underline{\lambda}}(q)=1.$$
Assume that the last part of the $a$-graph of $\underline{\lambda}$ 
is given by:
$$\underline{\nu}\underbrace{{\overset{i_s}{\underset{(r_1,p_1)}{\mapsto}}
...{\overset{i_s}{\underset{(r_{a_s},p_{a_s})}{\mapsto}}}}}
_{a_s}\underline{\lambda}.$$
Let $\underline{\mu}$ be a $d$-partition obtained  by removing
 $a_s$ $i_s$-nodes {}from $\underline{\lambda}$ and assume that
 $\underline{\mu}\neq{\underline{\nu}}$. Then, by construction
 of the $a-\sequence$ of residues and Proposition 
\ref{orderFLOTW}, the nodes $\gamma_l$ associated to 
$(r_l,p_l)$, $l=1,...,a_s$, are the lowest $i_s$-nodes of 
$\underline{\lambda}$ (with respect to the FLOTW order).
 If  $\underline{\mu}\neq{\underline{\nu}}$ then at least
 one of the lowest $i_s$-node is a node of $\underline{\mu}$. 
When this node is added, the $i_s$-node can be added to two 
parts, one to obtain $\underline{\mu}[t]$, the other to obtain 
 $\underline{\nu}[t]$ for some $t$. As we choose the higher 
node, this contradicts to the optimality of  $\underline{\nu}[t]$. 
Hence, we have to show that:
$$ \overline{N}_i^{b}(\underline{\nu},\underline{\lambda})=0.$$
There is no addable node of $\underline{\lambda}$ below the  $\gamma_l$
 and there is no  removable  node of $\underline{\nu}$ below the
  $\gamma_l$. Hence this  is clear.
\end{proof}

\begin{re}  This proposition shows that  $\underline{\lambda}$
 is a FLOTW $d$-partition if and only if there exists a sequence 
of residues  $\underbrace{i_1,...,i_1}_{a_1},\underbrace{i_{2},
...,i_{2}}_{a_{2}},...,\underbrace{i_s,...,i_s}_{a_s}$ such that 
$\underline{\lambda}$ is the minimal $d$-partition with respect 
to  $a$-function which appears in  the expression $\displaystyle
{f^{(a_s)}_{i_s}f^{(a_{s-1})}_{i_{s-1}} ...f^{(a_1)}_{i_1}
\underline{\emptyset}}$.

This result is similar to a conjecture of Dipper, James and 
Murphy (\cite{DJM}): in this paper, it was conjectured that
 $\underline{\lambda}$ is a Kleshchev $d${}-partition if and 
only if there exists a sequence of residues 
 $\underbrace{i_1,...,i_1}_{a_1},\underbrace{i_{2},...,i_{2}}_{a_{2}},
...,\underbrace{i_s,...,i_s}_{a_s}$ such that:
$$f^{(a_s)}_{i_s}f^{(a_{s-1})}_{i_{s-1}} ...f^{(a_1)}_{i_1}
\underline{\emptyset}=n_{\underline{\lambda}}\underline{\lambda}
+\sum_{\underline{\mu}\ntrianglerighteq{\underline{\lambda}}}
{n_{\underline{\mu}}\underline{\mu}},$$
 where $n_{\underline{\lambda}}\neq{0}$. The partial order
 $\trianglerighteq$ is called the dominance order and it is 
defined as follows.
$$\underline{\mu}\trianglerighteq \underline{\lambda} \iff 
\sum_{k=0}^{j-1} |\mu^{(k)}|+\sum_{p=1}^i \mu_p^{(j)}\geq  
\sum_{k=0}^{j-1} |\lambda^{(k)}|+\sum_{p=1}^i \lambda_p^{(j)}
\ \ \ (\forall j,i).$$
\end{re}

Now, the next proposition  follows from the same argument  
as \cite[Lemma 6.6]{LLT}, the proof  gives an explicit algorithm 
for computing the canonical basis, which is an analogue of the
 LLT algorithm (the details of the algorithm  will be published  elsewhere). 
\begin{prp}  Following the notations of paragraph \ref{combinatoire}
, the canonical basis elements of $\overline{\mathcal{M}}$ are of the
 following form:
$$\underline{\lambda}+\sum_{a(\underline{\mu})>a(\underline{\lambda})}
b_{\underline{\lambda},\underline{\mu}}(q)\underline{\mu},$$
  where $b_{\underline{\lambda},\underline{\mu}}(q)\in{q\mathbb{Z}[q]}$
 and  $\underline{\lambda}$ is a FLOTW $d$-partition.

Reciprocally, for any FLOTW $d$-partition $\underline{\lambda}$, 
there exists a canonical basis element of the above form.
\end{prp}
\begin{proof}
Let $\underline{\mu}\in{\Lambda^{1}_{\{e;v_0,...v_{d-1}\}}}$ be 
a FLOTW $d$-partition and assume that:
 $$a-\sequence (\underline{\lambda})=\underbrace{i_1,...,i_1}_
{a_1},\underbrace{i_{2},...,i_{2}}_{a_{2}},...,\underbrace
{i_s,...,i_s}_{a_s}.$$
Then, we define:
$$A(\underline{\lambda}):=f^{(a_s)}_{i_s}f^{(a_{s-1})}_{i_{s-1}}
 ...f^{(a_1)}_{i_1}\underline{\emptyset}.$$
By the previous proposition, we have:
$$A(\underline{\lambda})=\underline{\lambda}+ 
\sum_{a(\underline{\mu})>a(\underline{\lambda})}{c_{\underline{\lambda}
,\underline{\mu}}(q)\underline{\mu}},$$
where  $c_{\underline{\lambda},\underline{\mu}}(q)
\in{\mathbb{Z}[q,q^{-1}]}$.

The $A(\underline{\lambda})$ are linearly independant 
and so,  the following set is a basis of  the 
$\mathcal{U}_{\mathcal{A}}$-module  $\overline{\mathcal{M}}_{\mathcal{A}}$ 
generated by the empty $d$-partition with respect to the
 JMMO action:
$$\{A(\underline{\mu})\ |\ \underline{\mu}\in{\Lambda^{1}}\}.$$
Note that the vertices of the crystal graph of 
 $\overline{\mathcal{M}}$ are given by the FLOTW $l$-partitions.  
 Let $\{B_{\underline{\mu}}\ |\ \underline{\mu}\in{\Lambda^{1}}\}$
 be the canonical basis of  $\overline{\mathcal{M}}$. There exist 
Laurent polynomials $m_{\underline{\mu},\underline{\nu}}(q)$ 
such that, for all $\underline{\mu}\in{\Lambda^{1}}$, we \mbox{ have:}
$$B_{\underline{\mu}}=\sum_{\underline{\nu}\in{\Lambda^{1}}}
m_{\underline{\mu},\underline{\nu}}(q)A(\underline{\nu}).$$
Now, we consider the  bar involution on $\mathcal{U}_{\mathcal{A}}$, 
this is the $\mathbb{Z}$-linear automorphism defined by:
$$\overline{q}:=q^{-1},\qquad{\overline{k_h}=k_{-h},}
\qquad{\overline{e_i}=e_{i},}\qquad{\overline{f_i}=f_{i}.}$$
It can  be  extended  to  $\overline{\mathcal{M}}_{\mathcal{A}}$ 
by setting, for all $u\in{\mathcal{U}}_{\mathcal{A}}$:
$$\overline{u\underline{\emptyset}}:=\overline{u}\underline{\emptyset}.$$

 The $A(\underline{\mu})$ are clearly invariant under the bar 
involution and so are the  $B_{\underline{\mu}}$ (this is by 
the definition of the canonical basis, see \cite{Arilivre}). 
Therefore, we have:
$$m_{\underline{\mu},\underline{\nu}}(q)=
m_{\underline{\mu},\underline{\nu}}(q^{{}-1}),$$
for all $\underline{\mu}$ and $\underline{\nu}$ in $\Lambda^{1}$. 

Let $\underline{\alpha}$ be one of the minimal $d$-partitions 
with respect to the $a$-function such that 
$m_{\underline{\mu},\underline{\alpha}}(q)\neq{0}$. 
Then, the coefficient of $\underline{\alpha}$ in 
$B_{\underline{\mu}}$ is:
$$b_{\underline{\mu},\underline{\alpha}}(q):=
\sum_{\underline{\nu}\in{\Lambda^{1}}}m_{\underline{\mu},
\underline{\nu}}(q)c_{\underline{\nu},\underline{\alpha}}(q).$$
We have:
\begin{itemize}
\item $m_{\underline{\mu},\underline{\nu}}(q)=0$ if 
$a(\underline{\nu})<{a(\underline{\alpha})}$.
\item $c_{\underline{\nu},\underline{\alpha}}(q)=0$ 
if  $a(\underline{\nu})>{{a(\underline{\alpha})}}$ and 
if  $a(\underline{\nu})={{a(\underline{\alpha})}}$ and  
$c_{\underline{\nu},\underline{\alpha}}(q)\neq{0}$   
then $\underline{\nu}=\underline{\alpha}$.
\end{itemize}
Hence, we have:
$$b_{\underline{\mu},\underline{\alpha}}(q)=m_{\underline{\mu}
,\underline{\alpha}}(q)\neq{0}.$$
Now, as the FLOTW $d$-partitions are  labeling the vertices of
 the crystal graph of  $\overline{\mathcal{M}}_{\mathcal{A}}$ ,
 we have:
$$ B_{\underline{\mu}}=\underline{\mu}\ (\textrm{mod}\ q).$$
Therefore, if $\underline{\mu}\neq{\underline{\alpha}}$, we have:
$$b_{\underline{\mu},\underline{\alpha}}(q)\in{q\mathbb{Z}[q]}.$$
Moreover, we have:
$$m_{\underline{\mu},\underline{\alpha}}(q)=m_{\underline{\mu},
\underline{\alpha}}(q^{-1}).$$
Therefore, we have  $\underline{\mu}=\underline{\alpha}$ and 
$b_{\underline{\mu},\underline{\mu}}(q)=1$. To summarize, we obtain:
$$ B_{\underline{\mu}}=A(\underline{\mu})+\sum_{a(\underline{\nu})
>a(\underline{\mu})}m_{\underline{\mu},\underline{\nu}}(q)
A(\underline{\nu}).$$
The theorem follows now from Proposition \ref{asuite}.
\end{proof}

Using the correspondence between the canonical basis elements 
and the  indecomposable projective  modules, we can reformulate 
this result as follows.
\begin{thm}\label{THM} Let $\mathcal{H}_{R,n}:=\mathcal{H}_{R,n}
(v;u_0,...,u_{d-1})$ be the Ariki-Koike algebra of type
 $G(d,1,n)$ over $R:=\mathbb{Q}{(\eta_{de})}$ and assume  that
 the parameters are given by:
  $$v=\eta_e,\qquad{u_i =\eta_e^{v_j}\ \ j=0,...,d-1,}$$
where $\eta_e:=\operatorname{exp}(\frac{2i\pi}{e})$  and 
 $0\leq{v_0}\leq ...\leq v_{d-1}<e$.

Let $\{P_R^{\underline{\mu}}\ |\ \underline{\mu}\in{\Lambda^0}\}$ 
be the set of indecomposable  projective $\mathcal{H}_{R,n}$-modules.
 Then, there exists a bijection between the set of Kleshchev
 $d$-partitions and the set of FLOTW $d$-partitions:
$$j:\Lambda^0\to{\Lambda^1},$$
such that, for all  $\underline{\mu}\in{\Lambda^0}$:
$$e([P_R^{\underline{\mu}}]_p)=\[S^{j(\underline{\mu})}\]+
\sum_{\underline{\nu}\in{\Pi_n^{d}}\atop{a(\underline{\nu})>
a(j(\underline{\mu}))}}{d_{\underline{\nu},\underline{\mu}} 
 \[S^{\underline{\nu}}\]}.$$
\end{thm}
In particular, this theorem shows that the decomposition matrix 
 of $\mathcal{H}_{R,n}$ has a lower triangular shape with $1$ 
along the diagonal.

\begin{re}\label{afre} Keeping the notations of the above theorem, 
we can attach to each simple $\mathcal{H}_{R,n}$-module $M$ an 
$a$-value in the following way:
$$ \textrm{If } M=D_R^{\underline{\mu}} \textrm{ for } \underline{\mu}
\in{ \Lambda^0} \textrm{ then } a_M:=a_{S^{j(\underline{\mu})}}=
a(j(\underline{\mu})).$$
 Then, we have:
$$a_M=\min{\{a_V\ |\ d_{V,M}\neq{0}\}}.$$
We remark that this property was shown in \cite{geck3} for Hecke 
algebras where the $a$-function was defined in terms of 
Kazhdan-Lusztig basis.

  We have also the following property:  for all 
 $\underline{\nu}\in{\Lambda^1}$:
$$[S^{\underline{\nu}}]=[D^{j^{-1}(\underline{\nu})}]+
\sum_{\underline{\mu}\in{ \Lambda^0    }\atop{a(j(\underline{\mu})) 
<a(\underline{\nu})}}{d_{\underline{\nu},\underline{\mu}}
 [ D^{\underline{\mu}}]}. $$
\end{re}
In the following section, we see consequences of the above 
theorem on the representation theory of Hecke algebras of
 type $B_n$. 

\section{Determination of the canonical basic set for Hecke algebras
 of type $B_n$}\label{cons}

Let $W$ be a finite Weyl group, $S$ be the set of simple 
reflections and let $H$ be the associated Hecke algebra defined 
 over $\mathbb{Z}[y,y^{-1}]$ where $y$ is an indeterminate. 
Let $v=y^2$.        
Then, $H$ is defined by:
\begin{itemize}
\item generators: $\displaystyle{\left\{ T_{w}\ |\ w\in{W} \right\}}$,
\item relations: for  $s\in{S}$ and $w\in{W}$:
$$ T_s T_w=
    \begin{cases}
        T_{sw} & \textrm{if } l(sw)>l(w), \\
       vT_{sw}+(v-1)T_w & \textrm{if } l(sw)<l(w).\\\end{cases}$$
where  $l$ is the usual lenght function.
\end{itemize}
 Let $K=\mathbb{Q}(y)$ and   $\theta: \mathbb{Z}[y,y^{-1}] 
\rightarrow{k}$ be a specialization such that $k$ is the 
field of fractions of $\theta (\mathbb{Z}[y,y^{-1}])$. We
 assume that the characteristic of $k$ is either $0$ or a
 good prime for $W$.  Let $H_K:=K\otimes_{\mathbb{Z}[y,y^{-1}]} H$ 
and $H_k:=k\otimes_{\mathbb{Z}[y,y^{-1}]}  H$. Then, following 
\cite[Theorem 7.4.3]{geck1}, we have a well-defined decomposition
 map $d_{\theta}$ between the Grothendieck groups of $H_K$-modules 
and $H_k$-modules. For $V\in{\Irr{H_K}}$, we have :
$$d_{\theta}([V])=\sum_{M\in{\Irr{H_k}}}d_{V,M}[M]$$
where the numbers $d_{V,M}\in{\mathbb{N}}$ are the decomposition numbers
 (see \cite[Chapter 7]{geck1} for details).
                              
In \cite{geckrouq}, Geck and Rouquier have defined a 
``canonical basic set'' $\mathcal{B}\subset{\Irr{H_K}}$
 which leads to a natural parametrization of the set $\Irr{H_k}$.

\begin{thm}[Geck-Rouquier \cite{geckrouq}] We consider
 the following subset of $\Irr{H_K}$:
$$\mathcal{B}:=\{V\in{\Irr{H_K}}\ |\ d_{V,M}\neq{0}\ 
\textrm{and}\ a(V)=a(M)\ \textrm{for some}\ M\in{\Irr{H_k}}\}$$
Then, there exists a unique bijection $\mathcal{B}\leftrightarrow 
\Irr{H_k}$, $V\leftrightarrow \overline{V}$, such that the
 following two conditions hold:
\begin{itemize}
\item[(i)] For all $V\in{\mathcal{B}}$, we have 
$d_{V,\overline{V}}=1$ and $a(V)=a(\overline{V})$.
\item[(ii)] If $V\in{\Irr{H_K}}$ and $M\in{\Irr{H_k}}$ are 
such that $d_{V,M}\neq{0}$, then we have $a(M)\leq a(V)$, with equality only for $V\in{\mathcal{B}}$ and $M=\overline{V}$.
\end{itemize}
The set  $\mathcal{B}$ is called the canonical basic set.
\end{thm}

{}From now to the end of the paper, we assume that  $k$ 
is a field of characteristic $0$. Then,  $\mathcal{B}$ 
has been already determined for all specializations for 
 type $A_{n-1}$ in \cite{geck2},  for type $D_n$ and $e$ 
odd in \cite{geck3} and for type $D_n$ and $e$ even in 
\cite{a-fonction}. 

The aim of this section is to find $\mathcal{B}$ in the 
case of Hecke algebras of type $B_n$. Let $W$ be a Weyl 
group of type $B_n$ with the following diagram:
\\
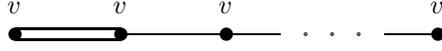
\begin{center}
\begin{picture}(240,20)
\put( 50,10){\circle*{5}}
\put( 50,8){\line(1,0){40}}
\put( 50,12){\line(1,0){40}}\put(207,18){$v$}
\put( 90,10){\circle*{5}}
\put( 90,10){\line(1,0){40}}\put( 47,18){$v$}
\put(130,10){\circle*{5}}
\put(130,10){\line(1,0){20}}
\put(160,10){\circle{1}}
\put(170,10){\circle{1}}\put( 87,18){$v$}
\put(180,10){\circle{1}}\put(127,18){$v$}
\put(190,10){\line(1,0){20}}
\put(210,10){\circle*{5}}
\end{picture}
\end{center}

Let $H$ be the correponding Hecke algebra over $\mathbb{Z}[y,y^{-1}]$.
 First, it is known that  $H_K$ is semi-simple unless $\theta (u)$ 
is a root of unity. In this case, we have $\mathcal{B}=\Irr{H_K}$.
 For $p\in{\mathbb{N}_{>0}}$  we put $\eta_p:=\textrm{exp}
{(\frac{2i\pi}{p})}$,
 then we can assume that $\theta (v)=\eta_e$ and that
 $k:=\mathbb{Q}(\eta_{2e})$.

The semi-simple algebra $H_K$ is an Ariki-Koike algebra 
with parameters $u_0=v$ and $u_1=-1$ over $K$.
Then, the simple  $H_K$-modules are given by the Specht 
modules $S^{\underline{\lambda}}$ which are labeled by 
the $2$-partitions $\underline{\lambda}$ of rank $n$. Let
   $\underline{\lambda}:=(\lambda^{(0)},\lambda^{(1)})$ 
be a $2$-partition and let $h^{(0)}$ (resp. $h^{(1)}$) be 
the height of $\lambda^{(0)}$ (resp.  $\lambda^{(1)}$). Let
 $r$ be a positive integer such that 
$r\geq{\textrm{max}\{h^{(0)},h^{(1)}\}}$.   Then, the $a$-value 
of the Specht module labeled by  $\underline{\lambda}:=
(\lambda^{(0)},\lambda^{(1)})$  is given by
$$\begin{array}{ll }a(\lambda^{(0)},\lambda^{(1)}):= 
&\displaystyle{ -\frac{1}{6}r(r-1)(2r+5)+\sum_{i=1}^r (i-1)
(\lambda^{(0)}_i+\lambda^{(1)}_i+1)}\\
&\displaystyle{ +\sum_{i,j=1}^r \textrm{min}{\{\lambda^{(0)}_i+1+r-i,
\lambda^{(1)}_j+r-j\}}},\end{array} $$
where we put  $\lambda^{(0)}_j:=0$ (resp.  $\lambda^{(1)}_j:=0$) 
if $h^{(0)}<j\leq r$ (resp. $h^{(1)}<j\leq r$). This formula is 
obtained by rewriting that in  Proposition \ref{af} with the above
  choice of parameters. Now, we have two cases to consider.\\
\\
a) Assume that $e$ is odd. In this case, we can apply  results of
 Dipper and James: they have shown that the simple  $H_k$-modules
  are given by the modules $D^{\underline{\lambda}}$ where
 $\underline{\lambda}=(\lambda^{(0)},\lambda^{(1)})$ is such that 
 $\lambda^{(0)}$ and  $\lambda^{(1)}$ are $e$-regular. 
A partition $\nu=(\nu_1,...,\nu_r)$ where $\nu_r>0$ is  $e$-regular
 if for all $i=1,...,r$, we can not have  $\nu_i=\nu_{i+1}=...=\nu_{i+e-1}$.
 For $\underline{\lambda}$  a $e$-regular $2$-partition 
 and $\underline{\mu}$ a $2$-partition, we denote by
 $d_{\underline{\mu},\underline{\lambda}}$ the corresponding
 decomposition number.

Moreover, Dipper and James have shown that the decomposition 
numbers of $H$ are determined by the decomposition numbers of
 a Hecke algebra of type $A_{n-1}$ in the following way 
(see \cite{DJ} and \cite{DipMat} for a more general case).

Let $0\leq{l}\leq{n}$ and let $H(\mathfrak{S}_l)$ be the generic
 Hecke algebra of type $A_{l-1}$, then $\theta$ determines a 
decomposition map between the Grothendieck groups of  
 $H_K(\mathfrak{S}_l )$ and $H_k(\mathfrak{S}_l )$. The simple
 modules of   $H_K(\mathfrak{S}_l)$  are given by Specht modules
 parametrized by partitions of rank $l$. The simple modules of  
 $H_k(\mathfrak{S}_l)$  are given by  the $D^{\lambda}$ labeled 
by the partitions $\lambda$ of rank $l$ which are $e$-regular. 
Let $d_{\lambda',\lambda}$ be the decomposition numbers of  
 $H(\mathfrak{S}_l)$  where  $\lambda$ runs over the set of 
$e$-regular partitions of rank $l$ and  $\lambda'$ over the set
 of  partitions of rank $l$. 

 Let $\lambda=(\lambda_1,...,\lambda_r)$ and $\mu=(\mu_1,...,\mu_r)$ 
be two partitions of rank $l$. Recall  that we write 
 $\lambda\trianglerighteq{\mu}$ if , for all $i=1,...,r$, we have:
$$\sum_{j=1}^i{\lambda_j}\geq{\sum_{j=1}^i{\mu_j}}.$$
Then, if $\underline{\lambda}=(\lambda^{(0)},\lambda^{(1)})$  
is a $e$-regular $2$-partition of rank $n$ and 
$\underline{\mu}=(\mu^{(0)},\mu^{(1)})$ a $2$-partition of
 rank $n$, we have, following \cite[Theorem 5.8]{DJ}:
\begin{equation*} 
 d_{\underline{\mu},\underline{\lambda}}=\begin{cases}
                        d_{\mu^{(0)},\lambda^{(0)}}
d_{\mu^{(1)},\lambda^{(1)}} & \text{if $|\mu^{(0)}|=|\lambda^{(0)}|$
 and $|\mu^{(1)}|=|\lambda^{(1)}|$},\\
0 & \text{otherwise.} \end{cases}
\end{equation*}
Assume that $|\mu^{(0)}|=|\lambda^{(0)}|$ and 
$|\mu^{(1)}|=|\lambda^{(1)}|$ and that  $d_{\mu^{(0)},\lambda^{(0)}}\neq{0}$,
  $d_{\mu^{(1)},\lambda^{(1)}}\neq{0}$, then following 
\cite[Theorem 3.43]{mathas} (result of Dipper and James),
 we have:
$$\mu^{(0)}\trianglelefteq{\lambda^{(0)}}\qquad{\textrm{and}}
\qquad{\mu^{(1)}\trianglelefteq{\lambda^{(1)}}}.$$
Then, as in \cite[Proposition 6.8]{geck3},   it is easy to see that:
$$a(\lambda^{(0)},\lambda^{(1)})\leq{a (\lambda^{(0)},\mu^{(1)})}
\leq{a (\mu^{(0)},\mu^{(1)})}.$$

For all $e$-regular $2$-partition $\underline{\lambda}$, we 
have $d_{\underline{\lambda},\underline{\lambda}}=1$.
Hence, it proves  the following proposition:
\begin{prp}  Keeping the above notations, assume that $W$ is 
a Weyl group of type $B_n$ and that  $e$ is odd, then the 
canonical basic set in bijection with $\operatorname{Irr}(H_k)$
 is the following one:
$$\mathcal{B}=\{S^{\underline{\lambda}}\ |\ \underline{\lambda}=
(\lambda^{(0)},\lambda^{(1)})\in{\Pi_n^2},\ \lambda^{(0)}\ 
\textrm{and}\ \lambda^{(1)}\ \text{are $e$-regular}\}.$$
\end{prp} 
b) Assume that $e$ is even. Then, by using the notations of
 paragraph \ref{not}, we put $v_0=1$ and $\displaystyle{v_1=\frac{e}{2}}$.
 Then, we have:
$$m^{(0)}=1\qquad{\textrm{and}}\qquad{m^{(1)}=0}.$$
Then, the Ariki-Koike algebra $\mathcal{H}_{K,n}$ over $K=\mathbb{Q}(y)$ have the following parameters:
$$u_0=y^2,$$
$$u_1=-1,$$
$$v=y^2.$$
This is nothing but the one parameter Hecke algebra $H_K$ 
of type $B_n$. If we specialize $y$ to $\eta_{2e}$, we obtain 
 the Hecke algebra $H_k$.

Thus, we are in the setting of Theorem  \ref{THM}. Thus, 
we obtain:
\begin{prp} Keeping the above notations, assume that $W$ is a 
Weyl group of type $B_n$ and that  $e$ is even, then the 
canonical basic set in bijection with $\Irr{H_k}$ is the following 
one:
$$\mathcal{B}=\{S^{\underline{\lambda}}\ |\ \underline{\lambda}=
(\lambda^{(0)},\lambda^{(1)})\in{\Lambda_{\{e;1,\frac{e}{2}\}}^1}\}.$$
We have $(\lambda^{(0)},\lambda^{(1)})\in{\Lambda^1_{\{e;1,\frac{e}{2}\}}}$
 if and only if:
\begin{enumerate}
\item we have:
\begin{align*}
&\lambda_i^{(0)}\geq{\lambda^{(1)}_{i-1+\frac{e}{2}}},\\
&\lambda^{(1)}_i\geq{\lambda^{(0)}_{i+1+\frac{e}{2}}}.
\end{align*}
\item  For all  $k>0$, among the residues appearing at the 
right ends of the length $k$ rows of   $\underline\lambda$, 
at least one element of  $\{0,1,...,e-1\}$ does not occur.
\end{enumerate}
\end{prp}

Thus, we obtain the parametrization of the canonical basic 
set for type $B_n$ in characteristic $0$ for all specializations. 
Note that  the canonical basic set  for the exceptional types 
can be easily deduced from the explicit tables of decomposition 
numbers obtained by Geck, Lux and M\"uller. Hence, the above  
results complete  the classification of the canonical basic set 
for all types and all specializations in characteristic $0$. \\

\textbf{Thanks.} The author wishes to thank Meinolf Geck for
 precious remarks and discussions. He also wants to thank 
Susumu Ariki, Bernard Leclerc and the referee for a lot  
of helpful comments and suggestions on this paper.

\begin{flushright}
\begin{tabular}{l}
 {\scshape Institut Girard Desargues,}\\
 {\scshape bat. Jean Braconnier, Universit\'e
Lyon 1,}\\
 {\scshape 21 av Claude Bernard,}\\
 {\scshape F--69622 Villeurbanne cedex, France.} \\
 $\ $\\
 Current Address:\\
 {\scshape Laboratoire de Math\'ematiques}\\
 {\scshape Nicolas Oresme,}\\
 {\scshape Universit\'e de Caen,}\\
 {\scshape BP 5186,}\\
 {\scshape F--14032 Caen cedex, France.} \\
E-mail: jacon@math.unicaen.fr
\end{tabular}
\end{flushright}

\end{document}